\documentclass[a4paper,11pt,reqno,twosides]{amsart}

\usepackage[backend=biber, giveninits, style=alphabetic, uniquelist=false, doi=false, url=false, isbn = false]{biblatex}
\renewcommand{\leq}{\leqslant}
\renewcommand{\geq}{\geqslant}
\usepackage{amssymb}

\newcommand{\dd}{\mathrm{d}}
\newcommand{\ii}{\mathrm{i}}

\newcommand{\Fc}{\mathcal F}

\newcommand{\Hc}{\mathcal H}

\newcommand{\Lc}{\mathcal L}

\newcommand{\Nc}{\mathcal N}

\newcommand{\Pc}{\mathcal P}
\newcommand{\Qc}{\mathcal Q}

\newcommand{\Sc}{\mathcal S}

\newcommand{\N}{\mathbb N}
\newcommand{\Nast}{\N^\ast} 

\newcommand{\Q}{\mathbb Q}
\newcommand{\R}{\mathbb R}

\newcommand{\Eb}{\mathbf{E}}

\newcommand{\Pb}{\mathbf{P}}
\newcommand{\Qb}{\mathbf{Q}}

\usepackage{mathrsfs}

\newcommand{\Bs}{\mathscr B}
\newcommand{\Cs}{\mathscr C}

\newcommand{\Xs}{\mathscr X}
\newcommand{\Zs}{\mathscr Z}

\newcommand{\Ff}{\mathfrak F}
\newcommand{\Hf}{\mathfrak H}
\newcommand{\Xf}{\mathfrak X}

\newcommand{\Tf}{\mathfrak T}
\newcommand{\Zf}{\mathfrak Z}

\usepackage{dsfont} \newcommand{\1}{\mathbf 1}

\renewcommand{\epsilon}{\varepsilon}
\renewcommand{\phi}{\varphi}

\newcommand{\KL}{\mathrm{KL}}
\newcommand{\TV}{\mathrm{TV}}
\newcommand{\Var}{\operatorname{Var}}

\newcommand{\defeq}{\vcentcolon=}
\newcommand{\eqdef}{=\vcentcolon}

\newcommand{\transposed}{\top}

\newcommand{\norm}[1]{\lVert #1 \rVert}

\usepackage{mathtools} \renewcommand{\complement}{\mathsf{c}}

\DeclareMathOperator{\sinc}{sinc}
\DeclareMathOperator{\si}{si}
\newcommand{\privpar}{\alpha}
\newcommand{\ppp}{\privpar^\prime}
\newcommand{\betaprime}{\beta^\prime}
\newcommand{\Laplace}{\Lc}

\newcommand{\fhat}{\widehat f}
\newcommand{\ftilde}{\widetilde f}
\newcommand{\ellhat}{\widehat \ell}

\newcommand{\Ktilde}{\widetilde K}

\newcommand{\Borel}{\Bs}

\newcommand{\hlow}{\underline h}
\newcommand{\hupp}{\overline h}
\newcommand{\hast}{h^\ast}
\newcommand{\hhat}{\widehat h}
\newcommand{\hstar}{h^\star}

\newcommand{\GP}{\mathrm{GP}}

\newcommand{\kappap}{\kappa^{\prime}}

\usepackage{lmodern}
\usepackage[x11names]{xcolor}
\usepackage[english]{babel}

\usepackage{csquotes}

\usepackage{amsmath}
\usepackage{amsfonts}
\usepackage{amssymb}
\usepackage{amsthm}
\usepackage{amscd}

\usepackage{bookmark} 

\usepackage{marvosym} 

\usepackage{mathtools} 

\usepackage[shortlabels]{enumitem}

\usepackage{stmaryrd} \newcommand{\llbr}{\llbracket}
\newcommand{\rrbr}{\rrbracket}

\usepackage{tabularx}
\usepackage[font=footnotesize]{caption}
\usepackage{colortbl}

\usepackage[notref, notcite, color, final]{showkeys}
\definecolor{refkey}{rgb}{0.0,0.0,0.0}
\definecolor{labelkey}{rgb}{0.0,0.0,0.0}

\theoremstyle{plain}
\newtheorem{theorem}{Theorem}[section]
\newtheorem{proposition}[theorem]{Proposition}
\newtheorem{lemma}[theorem]{Lemma}
\newtheorem{corollary}[theorem]{Corollary}

\theoremstyle{definition}
\newtheorem{definition}[theorem]{Definition}

\newtheorem{example}[theorem]{Example}

\theoremstyle{remark}
\newtheorem{remark}[theorem]{Remark}

\setlist{leftmargin=2em}

\newenvironment{smalign*}
 {\par$\!\aligned}
 {\endaligned$\par}

\allowdisplaybreaks 

\DeclareFieldFormat[article]{title}{#1} \DeclareFieldFormat[book]{title}{\itshape #1}
\DeclareFieldFormat[thesis]{title}{#1}
\DeclareFieldFormat[incollection]{title}{\itshape #1}
\DeclareFieldFormat[inproceedings]{title}{\itshape #1}
\DeclareFieldFormat[article]{author}{\color{blue} #1}
\DeclareFieldFormat[article]{pages}{#1.}
\DeclareFieldFormat[incollection]{pages}{#1.} 

\renewbibmacro{in:}{\ifentrytype{article}{}{\printtext{\bibstring{in}\intitlepunct}}}

\DeclareNameAlias{sortname}{last-first}
\DeclareNameAlias{default}{last-first}

\renewbibmacro*{volume+number+eid}{\printfield{volume}\setunit*{\addnbthinspace}\printfield{number}\setunit{\addcomma\space}\printfield{eid}}
\DeclareFieldFormat[article]{number}{\mkbibparens{#1}}

\renewbibmacro*{publisher+location+date}{\printlist{publisher}\iflistundef{location}
	{\setunit*{\addcomma\space}}
	{\setunit*{\addcomma\space}}\printlist{location}\setunit*{\addcomma\space}\usebibmacro{date}\newunit}

 \DeclareFieldFormat[article]{number}{}

\usepackage[paper=a4paper,left=30mm,right=30mm,top=30mm,bottom=30mm]{geometry}

\usepackage{hyperref}
\hypersetup{colorlinks=true,urlcolor=RoyalBlue4,linkcolor=RoyalBlue4,citecolor=RoyalBlue4,naturalnames=true,hypertexnames=false}

\addto{\captionsenglish}{}
\usepackage{amsaddr}

\title[Adaptive kernel density estimation under local privacy]{Pointwise adaptive kernel density estimation under local approximate differential privacy}
\author[M. Kroll]{Martin Kroll}
\address{CREST, ENSAE, Institut Polytechnique de Paris}
\email{prenom.nom[arobase]ensae.fr}

\thanks{The author gratefully acknowledges financial support from GENES and by the French National Research Agency (ANR) under the grant Labex Ecodec (ANR-11-LABEX-0047).}
\date{\today}
\subjclass[2010]{62G05 (primary), and 68P25 (secondary)}
\keywords{Kernel density estimation. Approximate local differential privacy. Reproducing kernel Hilbert space. Adaptive estimation. Lepski's method.}

\begin{document}

\begin{abstract}
We consider non-parametric density estimation in the framework of local approximate differential privacy.
In contrast to centralized privacy scenarios with a trusted curator, in the local setup anonymization must be guaranteed already on the individual data owners' side and therefore must precede any data mining tasks.
Thus, the published anonymized data should be compatible with as many statistical procedures as possible.
We suggest adding Laplace noise and Gaussian processes (both appropriately scaled) to kernel density estimators to obtain approximate differential private versions of the latter ones.
We obtain minimax type results over Sobolev classes indexed by a smoothness parameter $s>1/2$ for the mean squared error at a fixed point.
In particular, we show that taking the average of private kernel density estimators from $n$ different data owners attains the optimal rate of convergence if the bandwidth parameter is correctly specified.
Notably, the optimal convergence rate in terms of the sample size $n$ is $n^{-(2s-1)/(2s+1)}$ under local differential privacy and thus deteriorated to the rate $n^{-(2s-1)/(2s)}$ which holds without privacy restrictions.
Since the optimal choice of the bandwidth parameter depends on the smoothness $s$ and is thus not accessible in practice, adaptive methods for bandwidth selection are necessary and must, in the local privacy framework, be performed directly on the anonymized data.
We address this problem by means of a variant of Lepski's method tailored to the privacy setup and obtain general oracle inequalities for private kernel density estimators.
In the Sobolev case, the resulting adaptive estimator attains the optimal rate of convergence at least up to extra logarithmic factors.
\end{abstract}
 
\maketitle

\section{Introduction}

In the modern information era data are routinely collected in all areas of private and public life.
Although the availability of massive data sets is essential to answer important scientific and societal questions, the individual data owners (who may be individuals, households, research institutions, companies, \ldots) might refuse to share their, maybe sensitive, raw data with others.
Even more, in view of regularly reported data leaks, they may not even want to entrust their data to a central curator who stores the data and publishes anonymized summary statistics.
Finding ourselves in such a dilemma, the question whether and, if yes, how data analytics can still be performed is of special importance.
For the evaluation of this question, several aspects have to be taken into account.

Firstly, in absence of a trusted curator, privacy of the data has to be achieved already \emph{locally} at the individual data owners' level.
The $i$-th data holder takes its datum, say $X_i$, as the input of a privacy mechanism and creates an output $Z_i$ that is considered sufficiently anonymized, for instance, in the sense of any of the privacy definitions listed below.
For the purpose of the present paper, a privacy mechanism is a Markov kernel $Q_i$ between measurable spaces $(\Xf,\Xs)$ and $(\Zf,\Zs)$ generating $Z_i$ given $X_i=x$ according to the distribution $Q^{Z_i\mid X_i=x}$.
This definition of \emph{local} privacy is in contrast to the framework of \emph{centralized} or \emph{global} privacy where the trusted curator can take the whole data set $\{ X_1,\ldots,X_n \}$ to create an output $Z$.\footnote{Thus, the local privacy model can be seen as a proper submodel of the global one because the trusted curator can also mimic any conceivable procedure in the local model.}

Secondly, for the quantification of privacy, different solutions have been proposed so far
(see \cite{barber2014privacy}, Section~2 for a comprehensive overview of existing privacy definitions):
\begin{itemize}
 \item In this paper, we will exclusively work in the framework of $\privpar$-differential privacy and its generalization $(\privpar,\beta)$-differential privacy as defined in Definition~\ref{DEF:APPROX:DP} below.
 These two privacy definitions are also referred to as \emph{pure} and \emph{approximate differential privacy}, respectively.
 Originally, these concepts have been suggested for the anonymization of microdata tables in a global privacy setup, more precisely in a framework where queries are answered by a server that has direct access to the sensitive data~\cite{dwork2006differential,dwork2006calibrating,dwork2008differential}.
 In the statistics community, working under privacy constraints has been popularized in the past decade, amongst others, through the articles \cite{wasserman2010statistical,hall2013differential} (in the global setup) and \cite{duchi2018minimax} (in the local privacy setup).
 Another strict relaxation of pure differential privacy is \emph{random differential privacy} as introduced in \cite{hall2011random}.
 \item An alternative quantification of privacy can be given as follows:
 Let $\phi$ be a function from $[0,\infty]$ to $\R \cup \{ +\infty \}$ with $\phi(1)=0$.
 Then, the associated \emph{$\phi$-divergence} between two distributions $\Pb,\Qb$ is
 \[ D_\phi(\Pb |\!| \Qb) = \int \phi \left( \frac{\dd \Pb}{\dd \Qb} \right) \dd \Qb = \int \phi\left( \frac{p}{q} \right) q \dd \mu \]
 where $\mu$ is a measure such that $\Pb,\Qb \ll \mu$ and $p,q$ denote the corresponding Radon-Nikodym densities.
 Then, the mechanism $Q$ is called \emph{$\privpar$-$\phi$-divergence private} if
 \[ \sup_{x,x^\prime \in \Xf}  D_\phi(Q(\cdot|X=x) |\!| Q(\cdot|X=x^\prime)) \leq \alpha. \]
\end{itemize}
The intersection of these two concepts is non-empty: For instance, taking $\phi(x) = \lvert x-1\rvert/2$, the $\phi$-divergence $D_\phi(\Pb |\!| \Qb)$ is the total variation distance, and the resulting $\privpar$-$\phi$-divergence is equivalent to $(0,\beta)$-differential privacy.

\bigskip

\begin{figure}[t]
  \center
  \includegraphics[width=0.99\textwidth]{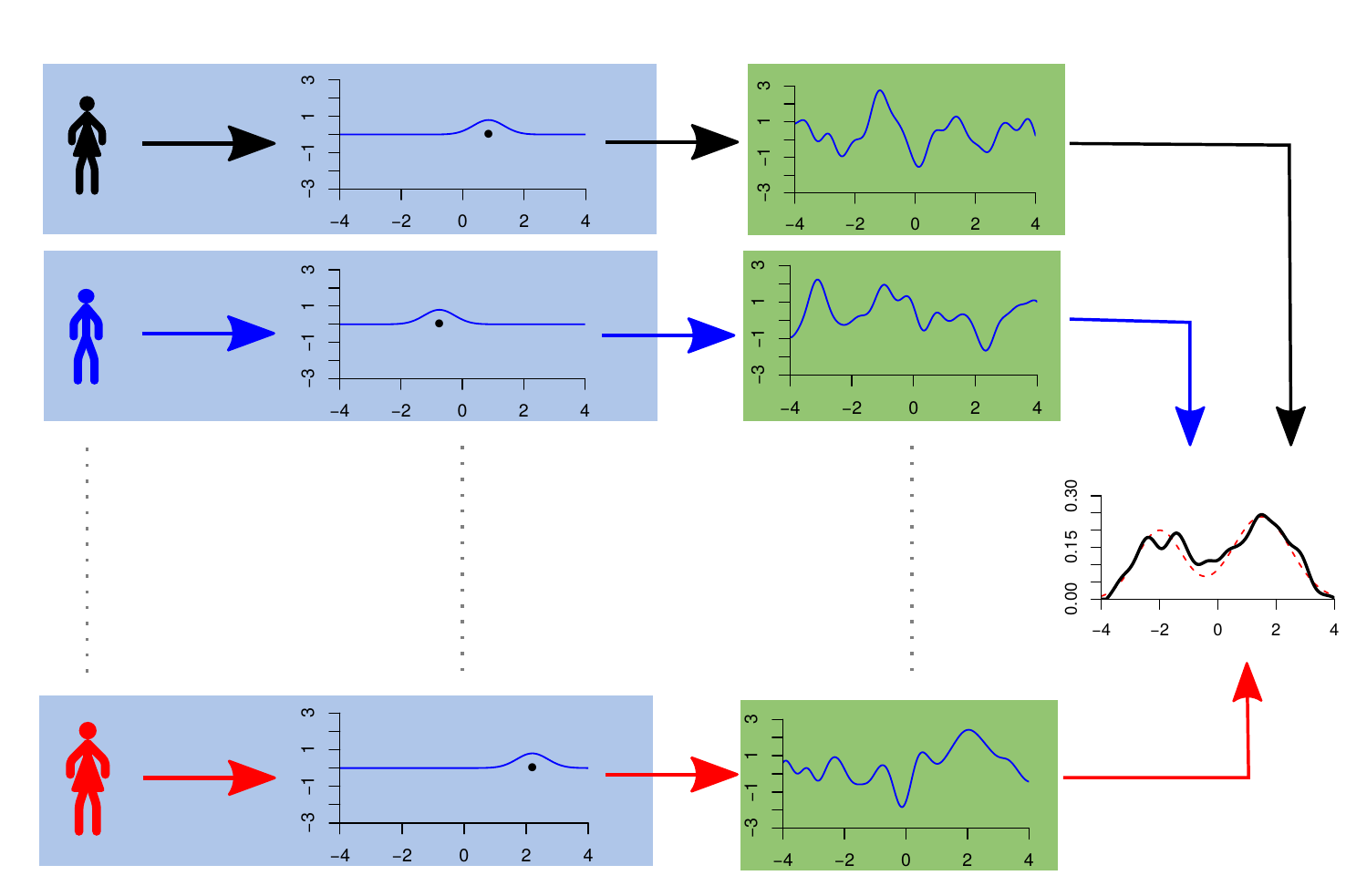}
  \captionsetup{width=\linewidth}
  \caption{General workflow of our procedure in the framework of univariate density estimation.
  Objects in blue boxes can only be observed by the respective data holder in the same box. Every data owner computes a kernel density estimator based only on its proper observation (given by the first coordinate of the black points). These estimators can be published after being perturbed by an appropriate centred Gaussian process (output in the pistachio coloured boxes).
  The pointwise mean of the private kernel density estimators (black solid curve) provides a natural estimator of the unknown density function (red dashed curve).}
  \label{FIG:WORKFLOW}
\end{figure}

Thirdly, the published data $Z_1,\ldots,Z_n$ should ideally be multi-purpose in the sense that they can serve as input data for several types of analyses.
Thus, when the unmasked data are for instance a sample from an unknown probability distribution, the anonymized data should contain as much information as possible about the whole distribution and not only about certain characteristics.
One main motivation for this work is to introduce novel methodology in the framework of density estimation that aims to address also this issue by proposing a local approximate differential private version of kernel density estimators, that is, the whole function $t \mapsto K((X_i-t)/h)/h$ for a bandwidth parameter $h > 0$ along with a study of their theoretical properties.
Figure~\ref{FIG:WORKFLOW} gives a foretaste and provides a graphical representation of the general workflow developed in this paper.

\subsection*{Roadmap of the article}
Throughout the article, we consider the paradigmatic example of non-parametric density estimation.
For the sake of simplicity, we assume that each of $n$ data holders $D_i$ observes a size-one sample $X_i$ from a (in this paper) univariate target density $f$, but refuses to share this observation.
In Section~\ref{SEC:PRIVACY}, we first introduce two mechanisms to estrange the value of an kernel density estimator at a single fixed point $t \in \R$.
The first approach is based on adding appropriately scaled Laplace noise.
The second approach is based on adding Gaussian noise and can be extended, using the ideas introduced in \cite{hall2013differential}, to an anonymized version of the whole kernel density estimator (as a function from $\R$ to $\R$) via perturbation by a suitable Gaussian process.
In Section~\ref{SEC:MINIMAX}, we consider estimation of the unknown density function under approximate differential privacy from a minimax point of view.
As the performance measure to evaluate arbitrary estimators, we consider the mean squared error at a fixed point.
Both the Laplace and the Gaussian perturbation approach attain the optimal rate $n^{-(2s-1)/(2s+1)}$ in terms of $n$ over Sobolev ellipsoids with smoothness index $s$ which is slower than the rate $n^{-(2s-1)/(2s)}$ in the setup without privacy constraints.
The Gaussian process approach, however, makes it possible to estimate the value of the density at any point of the observation window and not only at one single point that has to be chosen prior to the anonymization procedure.
In addition, this approach enables the statistician to perform any kind of analysis that plugs kernel density estimators into others estimators.
Investigating theoretical guarantees of such plug-in procedures, however, is outside the scope of this work and deferred to future research.

As usual for kernel density estimators, the choice of the bandwidth parameter is crucial.
In the considered minimax framework over Sobolev classes, the optimal order of the bandwidth that leads to a rate optimal estimator depends on the smoothness index $s$ which is typically unknown.
In Section~\ref{SEC:ADAPTATION}, we apply a Lepski scheme tailored to the privacy framework to overcome this problem and obtain an adaptive choice of the bandwidth.
This issue specifically arises in the local privacy setup since in the global framework the trusted curator can apply the existing plethora of methods for bandwidth selection on the unmasked data, and then only publish the resulting estimator with the adaptively determined bandwidth in its anonymized form.
In order to perform the Lepski scheme, any data owner has to publish the kernel density estimator not only for one single bandwidth but for a finite set of potential bandwidths.
Such a multiple output still guarantees the desired privacy condition provided that the additive noise is multiplied with a factor proportional to the number of potential bandwidths which is logarithmic in the number of data sources in our case.
We derive general oracle type inequalities for the estimator resulting from the Lepski procedure adapted to the privacy framework.
For the specific example of Sobolev ellipsoids, the rates of convergence are merely deteriorated by logarithmic factors with respect to the case of \emph{a priori} known smoothness.
 \section{Privacy mechanisms}\label{SEC:PRIVACY}

\subsection{Definition of approximate differential privacy}

Let $(\Xf,\Xs)$ and $(\Zf,\Zs)$ be measurable spaces.
A privacy mechanism is a Markov kernel $Q: \Xf \times \Zs \to [0,1]$ with the interpretation that, given original data $X=x$, an anonymized output is randomly drawn from the probability measure $Q(\cdot | X=x)$.
In the non-interactive setup that we are going to consider, we work under the following definition of \emph{approximate} or \emph{$(\privpar,\beta)$-differential privacy}.
\begin{definition}\label{DEF:APPROX:DP}
	Let $\privpar \geq 0, \beta \in [0,1]$.
	We say that $Z \sim Q(\cdot \mid X)$ is a local $(\privpar,\beta)$-differentially private view of $X$ if for all $x,x^\prime \in \Xf$, $A \in \Zs$ the estimate
	\begin{equation}\label{EQ:COND:APPROX:DP}
		Q(A \mid X = x) \leq \exp(\privpar) Q(A \mid X = x^\prime) + \beta,
	\end{equation}
	holds true.
\end{definition}
Let us emphasize that in Definition~\ref{DEF:APPROX:DP} the spaces $(\Xf,\Xs)$ and $(\Zf,\Zs)$ do not need to coincide.
In fact, in Example~\ref{EX:FUNC:GAUSSIAN} the space $(\Xf,\Xs)$ will be the real line equipped with its Borel sets and $(\Zf,\Zs)$ a measurable space of random functions.
In the literature, the case $\beta = 0$ is also referred to as $\privpar$-differential privacy or \emph{pure} differential privacy.
Evidently, the privacy condition \eqref{EQ:COND:APPROX:DP} becomes more restrictive for smaller values of the two parameter $\privpar$ and $\beta$.
Although Definition~\ref{DEF:APPROX:DP} smoothly bridges the cases $\beta = 0$ and $\beta > 0$, the classical anonymization techniques used for $\beta = 0$ and $\beta > 0$ are essentially different:
In the case $\beta = 0$, Laplace perturbation as well as randomization techniques as considered in \cite{duchi2018minimax,rohde2018geometrizing} can be used.
In the case $\beta > 0$, adding appropriately scaled Gaussian noise has been suggested in \cite{hall2013differential}.
However, as proved in \cite{holohan2015differential}, appropriately scaled Laplace noise can also lead to approximately differential private outputs (see Proposition~\ref{PROP:UNIV:LAPLACE:MECH} below).
In the sequel, we discuss how to achieve approximate differential privacy by means of these classical subroutines and how they can be extended to deal with functional data as well.

\subsection{Univariate output using Laplace noise}

First, we consider the case that both the input and the output of the privacy mechanism are univariate and real-valued, that is $(\Xf,\Xs)=(\Zf,\Zs)=(\R,\Borel(\R))$.
For this case, we consider Laplace perturbation which is also used to derive an upper bound in Section~\ref{SEC:MINIMAX}.
More precisely, let $Y_i=g(X_i) \in \R$ a quantity derived from the $X_i$ that should be masked.
Define the sensitivity of $g$ as
\[ \Delta(g) = \sup_{x,x^\prime \in \Xf} \lvert g(x) - g(x^\prime) \rvert. \]
Recall that the univariate Laplace distribution, denoted by $\Laplace(b)$, is given by the probability density function $p_b(x)=\frac{1}{2b} \exp(-\lvert x \rvert / b)$ (we include also the case $b=0$; then the Laplace distribution is, by convention, the Dirac measure concentrated at $0$).
In particular, the variance of an $\Laplace(b)$ distributed random variable is $2b^2$.
The following proposition establishes approximate differential privacy by Laplace perturbation.

\begin{proposition}[See~\cite{holohan2015differential}, Example~5]\label{PROP:UNIV:LAPLACE:MECH}
Let $\privpar > 0$, $\beta \in [0,1]$.
Then
\[ Z = g(X) + b \xi \]
with $\xi \sim \Laplace(1)$ for $b \geq \Delta(g)/(\privpar - \log(1- \beta))$
provides an $(\privpar,\beta)$-differential private view of $g(X)$ (and of $X$ as well).
\end{proposition}

A benefit of Proposition~\ref{PROP:UNIV:LAPLACE:MECH} in contrast to the often proposed perturbation by Gaussian noise to establish approximate differential privacy is that it allows to deal with the cases $\beta = 0$ and $\beta > 0$ by the same approach.
Moreover, letting the parameter $\beta$ vary permits natural interpretations:
If $\beta = 0$, the variance of $\sqrt 2 b \xi$ corresponds to the one that is usually encountered in the case of pure differential privacy.
When $\beta$ tends to one, the privacy constraint gets weaker and the variance of the centred noise $\sqrt 2 b \xi$ tends to $0$.
In the extreme case $\beta = 1$ it is even allowed to publish $g(X)$ directly.

We now introduce kernel density estimators as the guiding example that we have in mind for the function $g$ for the rest of the paper.

\begin{example}\label{EX:UNIV:LAPLACE}
Let $X_1,\ldots,X_n$ i.i.d.\,according to an unknown probability density function $f \colon \R \to \R$.
Let $t \in \R$ be fixed.
Then the $i$-th dataholder, who observes $X_i \in \R$, can compute
\[ K_h(X_i-t) \defeq \frac{1}{h} K \left( \frac{X_i - t}{h} \right) \]
for a bounded kernel function $K$, that is, $K \colon \R \to \R$ is integrable and $\int K(u) \dd u = 1$.
The quantity $K_h(X_i-t)$ will play the role of $g(X)$ in Proposition~\ref{PROP:UNIV:LAPLACE:MECH}.
By the triangle inequality $\Delta(K_h(\cdot-t)) \leq 2\lVert K \rVert_\infty/h$, and one can take any $b \geq 2 \lVert K \rVert_\infty/(h(\privpar - \log(1- \beta)))$ to obtain an approximate differential private view of $K_h(X_i-t)$.
Note that $t \in \R$ has been fixed in advance before the anonymization procedure.
\end{example}

\subsection{Multivariate output}

In principle, also multivariate output could be dealt with by adding independent Laplace noise to any of the components of the vector to be published.
In this case, both $\privpar$ and $\beta$ for each component have to be appropriately scaled in order to obtain the desired level of approximate differential privacy for the whole vector (the scaling can be carried out, for instance, as described in Lemma~\ref{LEM:COMP} below).
This approach, however, results in an increase concerning the Laplace noise added at any single point where the kernel density estimator is evaluated, and thus might deteriorate the performance of subsequent analyses more than necessary.
We do not further pursue this course here, since we will introduce a method for the anonymization of functional data that does not inflate the noise at single points in the next subsection.
Having stated this general method, we can, for instance, anonymize the whole function $\cdot \mapsto K_h(X_i-\cdot)$ in Example~\ref{EX:UNIV:LAPLACE}, and as a by-product we obtain $(\privpar,\beta)$-differential privacy for all pointwise evaluations $K_h(X_i-t)$, $t \in \R$ without any extra cost on the noise to be added.
To achieve anonymization of functional data, adding Gaussian processes with appropriately chosen covariance structure turns out to be convenient.
This idea has been originally suggested in \cite{hall2013differential}, but we state the essential steps here again for a clear exposition, and refer to \cite{hall2013differential} only for the proofs.
The first stopover on our way along the results from \cite{hall2013differential} is the following proposition that provides a condition under which approximate differential privacy of a vector is obtained by adding multivariate Gaussian noise with not necessarily uncorrelated components.
\begin{proposition}\label{PROP:PRIV:MULTI}
Let $\alpha > 0$, $\beta \in (0,1/2)$.
Let further $\Sigma \in \R^{m \times m}$ be a positive definite matrix and $g: \Xf \to \R^m$ for some $m \in \N^\ast$.
Assume that
\begin{equation}\label{EQ:COND:DELTA:MULTI}
   \sup_{x,x^\prime \in \Xf} \norm{ \Sigma^{-1/2} (g(x)-g(x^\prime)) }_2 \leq \Delta
\end{equation}
for all $x,x^\prime \in \Xf$.
Then, $Z$ defined via
\begin{equation*}
   Z = g(X) + \sigma \Xi, \qquad \Xi \sim \Nc(0,\Sigma),
\end{equation*}
is $(\privpar,\beta)$-differential private provided that
\begin{equation}\label{EQ:COND:SIGMA}
 \sigma \geq \frac{\Delta}{\privpar} \sqrt{2 \log(1/(2\beta)) + 2\privpar}.
\end{equation}
\end{proposition}

Proposition~\ref{PROP:PRIV:MULTI} will unfold its full potential in the next subsection where the condition \eqref{EQ:COND:DELTA:MULTI} will be reformulated appropriately.
For the univariate case (taking $m=1$), Proposition~\ref{PROP:PRIV:MULTI} directly provides a result similar to the one in Example~\ref{EX:UNIV:LAPLACE}, again with $t \in \R$ fixed before anonymization.

\begin{example}\label{EX:MULTI:GAUSSIAN}
We consider $K_h(X_i - t)$ as in Example~\ref{EX:UNIV:LAPLACE} and apply Proposition~\ref{PROP:PRIV:MULTI} for $m=1$ and $\Sigma=\begin{pmatrix}
	1
\end{pmatrix}$.
As in Example~\ref{EX:UNIV:LAPLACE},
\[ \sup_{x,x^\prime \in \R} \left\lvert \frac{1}{h} K \left( \frac{x-t}{h} \right) - \frac{1}{h} K \left( \frac{x^\prime-t}{h} \right) \right\rvert \leq \frac{2 \norm{K}_\infty}{h}, \]
and one can take $\Delta(K_h(\cdot - t)) = 2\norm{K}_\infty/h$ in \eqref{EQ:COND:DELTA:MULTI}.
Then, Proposition~\ref{PROP:PRIV:MULTI} guarantees that the $Z_i$, $i=1,\ldots,n$ defined through
\begin{equation*}
    Z_i = \frac{1}{h} K \left( \frac{X_i - t}{h} \right) + \frac{2  \norm{K}_\infty \sqrt{2\log (1/2\beta)+2\alpha}} {\privpar h} \, \xi_i, \qquad \xi_i \sim \Nc(0,1),
\end{equation*}
is an $(\privpar,\beta)$-differential private view for $\privpar, \beta > 1/2$ of $g_t(X_i)$ (and of $X_i$ as well).
\end{example}

\subsection{From multivariate to functional output}

The anonymization techniques used in Examples~\ref{EX:UNIV:LAPLACE} and \ref{EX:MULTI:GAUSSIAN} both suffer from the drawback that the output $Z_i$ provides information on the kernel density estimator $K_h(X_i-t)$ for one single $t$ only.
The aim of this section, based on Proposition~\ref{PROP:PRIV:MULTI} and ideas introduced in \cite{hall2013differential} in the context of global privacy, is to construct a privatized version of the whole function $t \mapsto K_h(X_i-t)$ by adding a suitable Gaussian process to the kernel density estimator.
As a consequence, the kernel density estimator anonymized in this vein can be evaluated at any single $t \in \R$.

For univariate and multivariate real-valued outputs of privacy mechanisms, the role of the $\sigma$-field $\Zs$ in Definition~\ref{DEF:APPROX:DP} is canonically taken by the Borel sets on $\R$ or $\R^m$.
In the case of functional output $Z \colon \Xf \to \R^m$ (where $\Xf$ is an arbritary set), its role is taken by the $\sigma$-field $\Cs$ which is generated by the cylinder sets
\[ C_{\Tf,B} = \{ f \colon \Xf \to \R : (f(t_1),\ldots,f(t_m)) \in B  \} \]
where $\Tf$ ranges over all finite sets $\Tf = \{ t_1,\ldots,t_m \} \subseteq \Xf$ and $B \in \Borel(\R^m)$.
The following result is a reformulation of Proposition~7 in \cite{hall2013differential} and we omit its proof.
See also Example~4 in \cite{holohan2015differential} for an alternative reasoning.
\begin{proposition}\label{PROP:PRIV:FUNC:GP}
Let $\Xi: \Xf \to \R$ be a sample path of a centred Gaussian process with covariance kernel $K : \Xf \times \Xf \to \R$.
For $t_1,\ldots,t_m \in \Xf$, consider the \emph{Gram matrix}
\begin{equation*}
 \Sigma_{t_1,\ldots,t_m} = \begin{pmatrix}
                             K(t_1,t_1) & \ldots & K(t_1,t_m)\\
                             \vdots & \ddots & \vdots\\
                             K(t_m,t_1) & \ldots & K(t_m,t_m)
                            \end{pmatrix}.
\end{equation*}
Let $X\colon \Xf \to \R$ be a (random) function in a function class $\Ff$.
Then, the release of
\begin{equation*}
  Z = X + \sigma \Xi
\end{equation*}
with $\sigma$ fulfilling \eqref{EQ:COND:SIGMA} is $(\privpar,\beta)$-differential private (with respect to $\Cs$) provided that
\begin{equation}\label{EQ:COND:DELTA:GP}
 \sup_{f,g \in \Ff} \sup_{m \in \Nast} \sup_{t_1,\ldots,t_m \in \Xf} \left\lVert \Sigma_{t_1,\ldots,t_m}^{-1/2} \begin{pmatrix}
                                                                                              f(t_1) - g(t_1)\\
                                                                                              \vdots \\
                                                                                              f(t_m) - g(t_m)
                                                                                             \end{pmatrix}
 \right\rVert_2 \leq \Delta
\end{equation}
where $\Delta$ is defined in \eqref{EQ:COND:DELTA:MULTI}.
\end{proposition}
The main question arising from Proposition~\ref{PROP:PRIV:FUNC:GP} is how the, on a first sight unhandy condition \eqref{EQ:COND:DELTA:GP}, might be verified.
The solution consists in transferring the problem into a reproducing kernel Hilbert space (RKHS) setup.
In fact, Proposition~\ref{PROP:PRIV:FUNC:GP} can be applied effectively when the random functions to be masked belong to the RKHS which corresponds to the covariance kernel of the Gaussian process $\Xi$.

In order to formulate this next result from \cite{hall2013differential}, we need to introduce some basic notation concerning the considered RKHS (we refer the reader to \cite{berlinet2004reproducing} for a comprehensive introduction to RKHS theory).
Let $K\colon \Xf \times \Xf \to \R$ be a positive definite kernel.
Recall that a real-valued kernel $K\colon \Xf \times \Xf \to \R$ is positive definite if
\begin{equation}\label{EQ:COND:PDK}
	\sum_{i,j=1}^k a_ia_j K(x_i,x_j) \geq 0
\end{equation}
holds for any $k \in \Nast$, $\{ a_1,\ldots,a_k \} \subseteq \R$, and $\{ x_1,\ldots,x_k\} \subseteq \Xf$.
For any $x\in \Xf$, define the function $K_x: \Xf \to \R$ by $K_x(\cdot) = K(x,\cdot)$.
Then the set
\begin{equation*}
 \Hf_0 \defeq \{ f : f = \sum_{i \in I} c_i K_{x_i} \text{ for some finite index set }I \}
\end{equation*}
is a pre-Hilbert space with respect to the norm $\norm{\cdot}_{\Hf}$ induced by the scalar product
\begin{equation*}
 \langle f,g \rangle_\Hf = \sum_{i\in I} \sum_{j\in J} c_i d_j K(x_i, y_j)
\end{equation*}
for $f=\sum_{i\in I} c_i K_{x_i}$, $g=\sum_{j\in J} d_j K_{y_j}$.
The RKHS corresponding to the kernel $K$ is the Hilbert space $\Hf$ resulting from the completion of $\Hf_0$ with respect to the RKHS norm $\norm{\cdot}_{\Hf}$.
The following two results are again taken from \cite{hall2013differential}.

\begin{proposition}[See \cite{hall2013differential}, Proposition~8]
 For $f \in \Hf$, where $\Hf$ is the RKHS corresponding to the kernel $K: \Xf \times \Xf \to \R$, and for any finite sequence $t_1,\ldots,t_m$ of distinct points from $\Xf$, we have
 \begin{equation*}
      \left\lVert \begin{pmatrix}
                  K(t_1,t_1) & \ldots & K(t_1,t_m)\\
                             \vdots & \ddots & \vdots\\
                             K(t_m,t_1) & \ldots & K(t_m,t_m)
                  \end{pmatrix}^{-1/2}
                  \begin{pmatrix}
                  f(t_1)\\
                  \vdots \\
                  f(t_m)
                  \end{pmatrix}
 \right\rVert_2 \leq \norm{f}_{\Hf}.
 \end{equation*}
\end{proposition}

\begin{corollary}[See \cite{hall2013differential}, Corollary~9]\label{COR:FUNC:OUTPUT:RKHS}
 For $X \in \Ff \subseteq \Hf$, the release of
 \[ Z = X + \sigma \Xi \]
 with $\sigma$ as in~\eqref{EQ:COND:SIGMA} is $(\privpar,\beta)$-differential private with respect to $\Cs$ provided that
  \begin{equation}\label{EQ:COND:DELTA:RKHS}
   \sup_{f,g \in \Ff} \norm{f-g}_\Hf \leq \Delta,
  \end{equation}
 and $\Xi$ is the sample path of centred Gaussian process with covariance kernel $K$ (given by the reproducing kernel of $\Hf$).
\end{corollary}
We now apply Corollary~\ref{COR:FUNC:OUTPUT:RKHS} to kernel density estimators.
\begin{example}\label{EX:FUNC:GAUSSIAN}
	In the case of univariate density estimation the $i$-th data holder observes $X_i$ drawn from the target density $f$, and we want him to be able to publish a approximately differential private version of the kernel density estimator
	\begin{equation*}
	   \ftilde_{i,h}(t) = \frac{1}{h} K \left( \frac{X_i-t}{h} \right), t \in \R,
	\end{equation*}
	based on his single observation $X_i$ only.
	In order to apply the above theory we have to assume that the kernel $K(x,y) = K(x-y)$\footnote{We slightly abuse notation by denoting both the kernel of the kernel density estimator and the corresponding kernel $\R \times \R \to \R$ given through $(x,y) \mapsto K(x-y)$ by the letter $K$.} is also a positive definite kernel.
	Under this additional assumption, Corollary~\ref{COR:FUNC:OUTPUT:RKHS} shows that the perturbed kernel density estimator
	\begin{equation*}
	   Z_{i,h}(\cdot) = \ftilde_{i,h}(\cdot) + \frac{\Delta}{\privpar} \sqrt{2 \log(1/(2\beta)) + 2\privpar} \Xi
	\end{equation*}
	where $\Xi$ a Gaussian process with kernel $hK_h(x,y)=K((x-y)/h)$ ensures $(\privpar,\beta)$-differential privacy provided that \eqref{EQ:COND:DELTA:RKHS} is satisfied.
	For instance, for the Gaussian kernel $K_{\text{Gauss}}(\cdot) = \exp(-(\cdot)^2/2h^2)$ we have
	\begin{align*}
	    \norm{K_h(x - \cdot) - K_h(x^\prime - \cdot)}_{\Hc}^2 = \frac{1}{2\pi h^2} (K_{\text{Gauss}}(0) + K_{\text{Gauss}}(0) - 2K_{\text{Gauss}}(x-x^\prime)) \leq \frac{1}{\pi h^2},
	\end{align*}
	and we can take $\Delta = 1/(\sqrt \pi h)$
	(the same argument working for any non-negative bounded kernel, and with a slight modification for any bounded kernel).
\end{example}
Let us emphasize that the property of positive definiteness is not satisfied for all kernels commonly used for kernel density estimators in non-parametric statistics.
In the following, we discuss some popular examples.
\begin{example}
	The \emph{rectangular kernel} given by
	\[ K_{\sqsubset\!\sqsupset}(x,y) \propto \1_{\{ \lvert x-y\rvert \leq 1 \} } \]
	for $x,y \in \R$ is \emph{not} positive definite.
	In order to see this, set $x_1=0$, $x_2=\frac{3}{4}$, $x_3=\frac{3}{2}$, $a_1=a_3=1$, and $a_2=-1$.
	Then
	\begin{align*}
		 \sum_{i=1}^3 \sum_{j=1}^3 a_i K_{\sqsubset\!\sqsupset}(x_i,x_j) a_j \propto 3-4 < 0,
	\end{align*}
	which contradicts the defining property \eqref{EQ:COND:PDK} of positive definite kernels.
\end{example}

\begin{example}
	The \emph{triangular kernel} given by
	\[ K_\triangle(x,y) \propto (1-\lvert x-y \rvert) \1_{  \{ \lvert x-y\rvert \leq 1 \} }  \]
	for $x,y \in \R$ is positive definite.
	This follows from the fact that kernels of the form
	\[ K(x,y) = \int_{\R^d} f(t+y)f(t+y)\dd t \]
	for $x,y \in \R^d$ with square integrable $f\colon \R^d \to \R$ are positive definite and
	\[ (1-\lvert x-y \rvert) \1_{ \{ \lvert x-y\rvert \leq 1 \} } = \int_\R \1_{[0,1/2]}(t+x) \1_{[0,1/2]}(t+y) \dd t. \]
 \end{example}

\begin{example}
	The \emph{Gaussian kernel}
  \[ K(x,y) \propto \exp(- \lvert x-y\rvert^2/2) \]
	and the \emph{exponential kernel}
	\[ K(x,y) \propto \exp(- \lvert x-y\rvert) \]
	are positive definite.
	These kernels of the form $K(x,y) \propto \exp(- |x-y|^\gamma)$ are positive definite if and only if $\gamma \in [0,2]$.
	This follows by combination of Theorem~2.2 and Exercise~2.13, (b) in \cite{berg1984harmonic}.
\end{example}

\begin{example}\label{EX:SINC}
	The \emph{$\sinc$ kernel} given by
  \[ K_{\sinc}(x,y) = \sinc(x-y) = \frac{\sin(\pi(x-y))}{\pi (x-y)} \]
  is positive semidefinite since the $\si$-function is the characteristic function of the uniform distribution on the interval $[-1,1]$.
  The $\sinc$-kernel attains also negative values but grant to the estimate $1 \geq \sinc(\cdot) \geq -0.3$ we have, in analogy to the calculation in Example~\ref{EX:FUNC:GAUSSIAN},
  \[ \lVert (K_{\sinc})_x - (K_{\sinc})_{x^\prime} \rVert^2_{\Hf} = \frac{1}{h^2} (K_{\sinc}(x,x) + K_{\sinc}(x^\prime,x^\prime) -2K_{\sinc}(x,x^\prime)) \leq \frac{2.6}{h^2} \]
	which yields a suitable bound for $\Delta$ in this example.
\end{example}

\begin{example}
	The \emph{Epanechnikov kernel}
	\[ K(x,y) = \frac{3}{4} (1-\lvert x-y\rvert^2) \1_{ \{ \lvert x-y\rvert \leq 1 \} } \]
	is \emph{not} positive definite. In order to see this, put $x_1=0$, $x_2=1/2$, $x_3=1$, $a_1=a_3=-0.9$ and $a_2=1$.
	Then,
	\begin{align*}
		\sum_{i=1}^3 \sum_{j=1}^3 a_i K(x_i,x_j) \overline a_j &= \frac{3}{4} \left[ 0.81 + 0.81 + 1 - 2 \cdot 0.9 \cdot 0.75 - 2 \cdot 0.9 \cdot 0.75 \right] = -0.08 < 0,
	\end{align*}
	in contradiction to the defining property \eqref{EQ:COND:PDK} of positive definite kernels.
\end{example}

\begin{example}
	The \emph{biweight kernel}
	\[ K(x,y) = \frac{15}{16} (1-\lvert x-y\rvert^2)^2 \1_{ \{ \lvert x-y\rvert \leq 1 \} } \]
	is \emph{not} positive definite.
  To see this, put $x_1=1/4$, $x_2 = -1/4$, $x_3 = -3/4$, and $x_4=1/2$.
	Then, consider the matrix $M=(K(x_i,x_j))_{i,j \in \llbr 1,4\rrbr}$.
	We have
	\[ \widetilde M \defeq 256 M = \begin{pmatrix}
		256 & 144 & 0 & 225\\
		144 & 256 & 144 & 49 \\
		0 & 144 & 256 & 0 \\
		225 & 49 & 0 & 256
	\end{pmatrix}, \]
	and the matrix $\widetilde M$ is not positive definite, since for $v = \begin{pmatrix}
			0.7 & -0.4 & 0.2 & -0.5
	\end{pmatrix}^\transposed$
	\[ v^\transposed \widetilde M v = -0.94 < 0. \]
\end{example}

\subsection{A composition lemma for approximate differential privacy}

For kernel density estimation, bandwidth selection is usually a delicate issue and so it is in our local privacy setup.
Whereas in the centralized setup existing methods can be applied by the trusted curator on the unmasked data, this is not possible in our local setup.
Thus the data holders have to publish versions of the kernel density estimator for different bandwidths, and one has to adapt general strategies from the non-private framework to the one with local approximately differential private data.
To do this under our privacy constraint it is necessary to understand how multiple outputs influence the defining condition of approximate differential privacy.
The following lemma provides a result of this flavour and is known in the research literature on privacy for statistical databases.
The setup is the following:
Given the unmasked datum $X$, the data owner does not only want to publish $Z_1=Z_1(X)$ but also $Z_2=Z_2(X)$, i.e., the vector $(Z_1,Z_2)$.
The following result tells us how $\privpar$ and $\beta$ for the single components have to be scaled in order to obtain $(\privpar, \beta)$-differential privacy for multiple outputs.
\begin{lemma}[Composition lemma for $(\privpar,\beta)$-differential privacy]\label{LEM:COMP}
  Let $Z_i$, $i=1,2$ be $(\privpar_i,\beta_i)$-differential private and conditionally (on $X$) independent views of $X$, respectively.
	Then $Z=(Z_1,Z_2)$ is an $(\privpar_1 + \privpar_2,\beta_1 + \beta_2)$-differential private view of $X$.
\end{lemma}

Of course, Lemma~\ref{LEM:COMP} can be successively applied.
For instance, if we want to publish $Z_{i,h}$ from the above examples for different $h$ in a finite set $\Hc$, then $\privpar$ and $\beta$ should be replaced with $\ppp = \privpar/\# \Hc$ and $\betaprime = \beta/\# \Hc$, respectively, in order to get differential privacy for $Z = (Z_{i,h})_{h \in \Hc}$.
 \section{Private minimax estimation}\label{SEC:MINIMAX}

Minimax theory provides a standard framework to study convergence properties of estimators in non-parametric statistics~\cite{tsybakov2009introduction}.
In this section, we apply this general toolbox to the specific case of density estimation under privacy constraints.
For fixed $t \in \R$ and any estimator $\ellhat$ of the linear functional $f(t)$ based on the private views $Z=\{ Z_1,\ldots,Z_n \}$, we study its mean squared error
\begin{equation*}
    \Eb[(\ellhat - f(t))^2].
\end{equation*}
The guiding  principle of minimax theory is to look for estimators that perform best in a worst-case scenario.
However, due to the privacy framework, we have not only the freedom of choosing the estimator $\ellhat$ but also the privacy mechanism $Q$ that generates the private outputs.
Hence, following \cite{duchi2018minimax}, classical minimax theory has to be adapted and a natural quantity to consider is the private minimax risk
\begin{equation*}
 \inf_{\substack{\ellhat \in \sigma(Z)\\Q \in \Qc_{\privpar,\beta}}} \sup_{f \in \Pc} \, \Eb[(\ellhat - f(t))^2]
\end{equation*}
where $\Pc$ is some function class containing probability densities and the infimum is taken over all local $(\privpar,\beta)$-differential private Markov kernels $Q \in \Qc_{\privpar,\beta}$ and all estimators based on the local approximate differential private views $Z$ of the corresponding original sample $X_1,\ldots,X_n$.
We specify the function class $\Pc$ by so called Sobolev ellipsoids $\Sc(s,L)$ that we define for $s > 1/2$ and $L>0$ by means of
\begin{equation*}
  \Sc(s,L) = \{ f \colon \R \to [0,\infty) : \int f(x) \dd x = 1, \int \lvert \Fc[f](\omega) \rvert^2 \lvert \omega \rvert^{2s} \dd \omega \leq 2 \pi L^2 \},
\end{equation*}
which, for $s \in \Nast$, is equivalent to the definition
\begin{equation*}
  \Sc(s,L) = \{ f \colon \R \to [0,\infty) : \int f(x)\dd x = 1, \int (f^{(s)}(x))^2 \dd x \leq L^2 \}.
\end{equation*}
In the first definition, $\Fc[f]$ denotes the Fourier transform of the density $f$, in the second one $f^{(s)}$ denotes the weak $s$-th derivative of $f$.

\subsection{Upper bound}

We first derive an upper bound on the minimax risk by specializing both the privacy mechanism $Q \in \Qc_{\privpar, \beta}$ and the estimator of $f(t)$.
Concerning the privacy mechanism, we consider the mechanisms mapping $X_i$ to private views $Z_{i,h}$ of $K_h(X_i - t)$ from Section~\ref{SEC:PRIVACY} for one single $h>0$.
More precisely, we consider the Laplace mechanism given through
\begin{equation}\label{EQ:OBS:LAPLACE}
  Z_{i,h}(t) = K_h(X_i-t) + \underbrace{\frac{2\lVert K\rVert_\infty}{h(\privpar - \log(1-\beta))}}_{\eqdef C^{\Laplace}_{\privpar \beta}/(\sqrt 2 h)}\, \xi_{i,h}, \qquad \xi_{i,h} \text{ i.i.d.} \sim \Laplace(1),
\end{equation}
and the Gaussian process mechanism given through
\begin{equation}\label{EQ:OBS:GP}
  Z_{i,h}(t) = K_h(X_i-t) + \underbrace{\frac{\Delta^\prime \sqrt{2\log(1/(2\beta)) + 2\privpar}}{h\privpar}}_{\eqdef C^{\GP}_{\privpar, \beta}/h}\, \Xi_{i,h}
\end{equation}
where $\Xi_{i,h}$ are i.i.d.\,Gaussian processes with covariance kernel $K((x-y)/h)$ and $\Delta^\prime$ is an upper bound on $\lVert (hK_h)_x - (hK_h)_{x^\prime} \rVert_\Hf$ for $x,x^\prime \in \R$.
Given $Z_{1,h},\ldots,Z_{n,h}$ as in \eqref{EQ:OBS:LAPLACE} or \eqref{EQ:OBS:GP}, a natural estimator of $f(t)$ is given by
\begin{equation}\label{EQ:DEF:FHAT}
	\fhat_h(t) = \frac{1}{n} \sum_{i=1}^{n} Z_{i,h}(t).
\end{equation}
The following proposition provides an upper risk bound for this estimator specialized with the $\sinc$-kernel over the Sobolev ellipsoids $\Sc(s,L)$ introduced above.

\begin{proposition}\label{PROP:UPPER}
 Consider the kernel density estimator $\fhat_h(t)$ for some fixed $t \in \R$ where the kernel used in the anonymization procedure \eqref{EQ:OBS:LAPLACE} or \eqref{EQ:OBS:GP} is the $\sinc$-kernel from Example~\ref{EX:SINC}.
 Then, for any $s > 1/2$,
 \[ \sup_{f \in \Sc(s,L)} \Eb [(\fhat_h(t) - f(t))^2] \leq C \left[ h^{2s-1} + \frac{1}{nh} + \frac{1}{nh^2} \right] \]
 for some $C=C(\privpar,\beta,L,s,\lVert f \rVert_\infty,K_{\sinc})$.
 In particular, setting $h=h^\star$ with $h^\star \asymp n^{-1/(2s+1)}$, we obtain
 \[ \sup_{f \in \Sc(s,L)} \Eb [(\fhat_{h^\star}(t) - f(t))^2] \lesssim n^{-\frac{2s-1}{2s+1}}. \]
\end{proposition}

Since the noise added by the privacy mechanisms is centred, the bias term in the proof of Proposition~\ref{PROP:UPPER} remains unchanged in comparison to the standard setup without privacy constraints.
However, the variance term changes due to the additional Laplace or Gaussian noise, respectively, and the classical variance term $1/(nh)$ is joined by the additional term $1/(nh^2)$ which is of higher order for $h \to 0$.
Consequently, the optimal bandwidth is no longer of order $n^{-1/(2s)}$ as in the standard setup but of the larger order $n^{-1/(2s+1)}$.
However, consistency of $\fhat_h$ is already guaranteed if $h \to 0$ and $nh^2 \to \infty$ simultaneously (in the standard density estimation setup one only needs $nh \to \infty$ in addition to $h \to 0$).

\subsection{Lower bound}

The following result states a lower bound over Sobolev ellipsoids in the case of pure differential privacy ($\beta = 0$).
\begin{proposition}\label{PROP:LOWER}
  Let $\privpar > 0$ arbitrary.
  Then,
	\begin{equation*}
		\inf_{\substack{\ellhat \in \sigma(Z) \\ Q \in \Qc_{\alpha,0}}} \sup_{f \in \Sc(s,L)} \Eb [(\ellhat - f(t))^2] \geq C(\privpar) n^{-\frac{2s-1}{2s+1}}
	\end{equation*}
	where $C(\privpar)> 0$ depends on the privacy parameter, and the infimum is taken over all estimators $\ellhat$ based on private views $Z_1,\ldots,Z_{n}$ and privacy mechanisms providing $(\privpar,0)$-differential privacy.
\end{proposition}

\begin{remark}
  The lower bound of Proposition~\ref{PROP:LOWER} still holds true when one allows a slight amount of interaction between the data holders, namely when the distribution of every $Z_i$ is determined by $X_i$ and the previously masked values $Z_1,\ldots,Z_{i-1}$.
  The proof remains the same because the data processing inequality (14) from \cite{duchi2018minimax} still holds true in this more general setup.
\end{remark}

Proposition~\ref{PROP:LOWER} shows that, regarding the privacy parameter $\privpar$ as an \emph{a priori} fixed constant, the estimators $\fhat_h(t)$ from Proposition~\ref{PROP:UPPER} attain the optimal rate $n^{-(2s-1)/(2s+1)}$ in terms of $n$ under pure local differential privacy.
Recall that without privacy restrictions the optimal rate over Sobolev ellipsoids is $n^{-(2s-1)/(2s)}$ (as mentioned in \cite{butucea2001exact}, this rate can, other than by a reduction scheme as used in our proof, be easily obtained via the theory developed in \cite{donoho1992renormalization}, see also \cite{tsybakov1998pointwise}).
In this work, we consider the parameters $\privpar$, $\beta$ as fixed and are interested in the behaviour of the rate as a function of $n$ only but remarks concerning $\privpar$ analogous to the ones made in \cite{butucea2019local} could be made (as in that paper, $\privpar$ and $\beta$ could also be allowed to vary with $n$).
The optimal behaviour, however, of the rates in terms of the privacy parameters $\privpar$ and $\beta$, especially if $\beta > 0$, remains an open issue.
 \section{Adaptation to unknown smoothness}\label{SEC:ADAPTATION}

The estimators of the previous section are not completely satisfying since the optimal choice $\hstar_n$ of the bandwidth, as usually in non-parametric statistics, depends on \emph{a priori} knowledge of the smoothness of the unknown function $f$.
Such knowledge is usually not available in practise.
At least, using the Gaussian process perturbation approach we relieved ourselves from the drawback of the Laplace method that one can privatize only one functional of the form $f(t)$ for one single $t$ that has to be fixed even before the anonymization.
Note that this drawback is, for instance, also present in the mechanisms suggested in \cite{rohde2018geometrizing}.
From this point of view, anonymization of the whole kernel density estimator via this approach should be preferred.

The purpose of this section is to address the remaining issue of adapting to the unknown smoothness of $f$.
In order to tackle this problem, we use a variant of Lepski's method (see~\cite{lepski1997optimal} for a general account in the Gaussian white noise model, and \cite{cavalier2001tomography} for an application to a tomography problem whose concise presentation has inspired our one).
Recall again that the necessity of novel methodology for adaptive estimation is specific for the setup of local privacy since in the global case the trusted curator can choose the bandwidth in an adaptive way using all the data $X_1,\ldots,X_n$ and, as a consequence, can build on the existing plethora of methods and theoretical results for this standard case; hence bandwidth selection does not provide any additional difficulty for centralized privacy since only the final output is anonymized.
In our local setup, where the data owners publish their data prior to any data analysis, adaptation must be addressed separately.
Note that the problem of adaptation has, to the best of the author's knowledge, only been addressed in the recent work~\cite{butucea2019local} so far, where the authors use wavelet estimators for density estimation on a compact interval.
The approach in that paper is thus conceptionally different from the one presented in the sequel.

We will apply Lepski's method both on observations \eqref{EQ:OBS:LAPLACE} where $t \in \R$ has been fixed \emph{a priori} and on pathwise observations \eqref{EQ:OBS:GP} from the Gaussian process approach that we evaluate at the point $t \in \R$ of interest.
In order to apply Lepski's method, the observations \eqref{EQ:OBS:LAPLACE} and \eqref{EQ:OBS:GP} must be available for different values of the bandwidth parameter $h$, say $h \in \Hc_n$.
This can be realized using Lemma~\ref{LEM:COMP} provided that the privacy parameters $\privpar$ and $\beta$ are appropriately scaled.
Thus, we can assume that $Z_{i,h}(t)$ are accessible for any $i \in \llbr 1,n\rrbr$ and $h \in \Hc_n$ if we replace $\privpar$ and $\beta$ by $\ppp=\privpar/\# \Hc_n$ and $\betaprime = \beta/\# \Hc_n$, respectively.
For any $h \in \Hc_n$ and $t \in \R$, we can then consider the estimator defined in \eqref{EQ:DEF:FHAT}.
In our case, we define the set of potential bandwidths by a geometrid grid,
\begin{equation*}
 \Hc_n = \{ h \in [\hlow_n, \hupp_n] : h = a^{-j}\hupp_n, j \in \N \},
\end{equation*}
where $a>1$ is a fixed constant, $\hupp_n$ is such that $a \log(\hupp_n \sqrt n)/\sqrt n \leq \hupp_n \leq 1$, and $\hlow_n$ satisfies $\hlow_n = (\log(\hupp_n \sqrt n) \vee 1)/\sqrt n$.
For $h \in \Hc_n$ and some $M > 0$, define\footnote{In the sequel, we write $C_{\ppp \betaprime}$ for both $C_{\ppp \betaprime}^\Laplace$ and $C_{\ppp \betaprime}^{\GP}$.}
\begin{equation*}
v^2(h) = \frac{M \int K^2(u)\dd u}{nh} + \frac{C_{\ppp\betaprime}^2}{nh^2}
\end{equation*}
where $C_{\ppp\betaprime}$ is defined as in Section~\ref{SEC:MINIMAX}.
The proof of \ref{PROP:UPPER} shows that
\[ \Var (\fhat_h) \leq v^2(h) \]
if $\norm{f}_\infty \leq M$.
Put $\lambda(h) = \max ( 1, ( \kappa \log(\hupp_n/h) )^{1/2} )$ with $\kappa$ being a sufficiently large constant (an explicit value can be determined from the proof of Theorem~\ref{THM:ADAPTATION}) and define
\begin{equation}\label{EQ:DEF:HAST}
 h^\ast_n = h_n^\ast(t,f) = \max \{ h \in \Hc_n : \lvert f_\eta(t) - f(t) \rvert \leq \frac{v(h)\lambda(h)}{2} \text{ for all }\eta \in \Hc_n, \eta \leq h \}.
\end{equation}
If the set in the definition of $\hast_n$ is empty, we set $\hast_n = \hlow_n$ by convention.
However, in the proof of Proposition~\ref{PROP:UPPER:ORACLE} we will show that this set is non-empty for $n$ large enough.
The bandwidth $h_n^\ast$ is an oracle in the sense that it is not accessible by the statistician since it depends on the unknown parameter $f$.
The definition of $\hast_n$ provides some kind of ideal criterion:
The bandwidth $h$ is increased along the grid $\Hc_n$ as long as the bias term $\lvert f_\eta(t) - f(t) \rvert$ it is bounded by the 'rate' $v(h)\lambda(h)$, a procedure that aims at mimicking the classical bias-variance tradeoff.
In order to state a risk bound for the pseudo estimator $\fhat_{\hast_n}$, we further define
\begin{equation*}
 r_n(t,f) = \inf_{\hlow_n \leq h \leq 1} \left[ \sup_{0 \leq \eta \leq h} (f_\eta(t) - f(t))^2 + \frac{M \int K^2(u)\dd u\log(n)}{nh} + \frac{C_{\ppp\betaprime}^2 \log(n)}{nh^2} \right].
\end{equation*}

\begin{proposition}\label{PROP:UPPER:ORACLE}
Consider the pseudo-estimator $\widehat f_{h^\ast_n}$ defined via~\eqref{EQ:DEF:FHAT} and \eqref{EQ:DEF:HAST} where $\privpar$ and $\beta$ are replaced with $\ppp$ and $\betaprime$, respectively.
Assume that
\begin{align}\label{EQ:ASS:BOCHNER}
  \lim_{h \to 0} \frac{1}{h} \int K \left( \frac{x-t}{h} \right) f(x) \dd x = f(t).
\end{align}
Consider $\hupp_n = 1$. Then, for $n$ sufficiently large,
\[ \Eb [ ( \fhat_{\hast_n}-f(t) )^2 ] \leq \frac{5}{4}v^2(\hast_n)\lambda^2(\hast_n) \leq C(a) r_n(t,f) \]
uniformly for all $f$ with $\lVert f\rVert_\infty \leq M$.
\end{proposition}

\begin{remark}
  Assumption~\eqref{EQ:ASS:BOCHNER} is satisfied in many cases.
  For instance, if $\int \lvert K(u) \rvert \dd u < \infty$, then \eqref{EQ:ASS:BOCHNER} is a special case of Bochner's lemma (see~\cite{tsybakov2004introduction}, Lemma~1.1).
  However, the $\sinc$-kernel is not absolutely integrable and thus Bochner's lemma cannot be applied.
  In this case, one can alternatively assume that $f$ belongs at least to some Sobolev space $\Sc(s,L)$ for some $s > 1/2$.
  Then, the analysis of the bias term as in the proof of Proposition~\ref{PROP:UPPER} guarantees the validity of \eqref{EQ:ASS:BOCHNER}.
\end{remark}

The pseudo estimator $\fhat_{\hast_n}$ is a stopover on our road to an adaptive estimator.
We now construct a genuine estimator of $f$ that aims at mimicking this oracle.
For this, we first define
\begin{align*}
  v^2(h,\eta) &= \frac{M}{n}\int (K_h(u) - K_\eta(u))^2 \dd u + \frac{C_{\ppp\betaprime}^2}{nh^2} + \frac{C_{\ppp\betaprime}^2}{n\eta^2}.
\end{align*}
Then, calculations similar to those in the proof of Proposition~\ref{PROP:UPPER} show that
\[ \Var(\fhat_h - \fhat_\eta) \leq v^2(h,\eta) \]
if $\norm{f}_\infty \leq M$.
For $h, \eta \in \Hc_n$, put
\[ \psi(h,\eta) = v(h)\lambda(h) + v(h,\eta) \lambda(\eta). \]
Then, we define an adaptive choice of the bandwidth parameter by
\begin{equation}\label{EQ:DEF:HHAT}
 \widehat h_n = \max \{ h \in \Hc_n : \lvert \widehat f_h(t) - \widehat f_\eta(t) \rvert \leq \psi(h,\eta) \text{ for all }\eta \leq h, h \in \Hc_n \}.
\end{equation}

This choice of the bandwidth is well-defined since the maximum is taken over a non-empty set.
The definition of $\hhat_n$ is characteristic for Lepski's method \cite{lepski1990problem}, and the motivation of this procedure is neatly described in \cite{cavalier2001tomography}, p.~67:
One chooses the largest bandwidth $h$ such that the difference between the two estimators $\fhat_h$ and $\fhat_\eta$ is not too large (in the sense of \eqref{EQ:DEF:HHAT}) for all $\eta \leq h$.
Evidently, the motivation of this procedure is to mimick the trade-off between squared bias and variance in a purely data-driven manner.
Note also that \eqref{EQ:DEF:HHAT} provides, as well as the oracle version~\eqref{EQ:DEF:HAST}, a local choice of the bandwidth in the sense that $\hhat_n$ depends on $t$.
Such a local criterion might result in a better adaptation to spatial inhomogeneity of the target density than global selection rules.

\begin{theorem}\label{THM:ADAPTATION}
 Consider the estimator $\fhat_{\hhat_n}$ defined via~\eqref{EQ:DEF:FHAT} and \eqref{EQ:DEF:HHAT} where $Z_{i,h}(t)$ for $h \in \Hc_n$ are defined via~\eqref{EQ:OBS:LAPLACE} or \eqref{EQ:OBS:GP} with $\privpar$ and $\beta$ replaced with $\ppp$ and $\betaprime$, respectively.
 Then, uniformly for all $f$ with $\lVert f\rVert_\infty \leq M$,
 \[ \Eb[(\fhat_{\hhat_n}(t) - f(t))^2] \leq C(a)v^2(\hast_n)\lambda^2(\hast_n). \]
 As a consequence, taking $\hupp_n = 1$, we obtain
 \[  \Eb[(\fhat_{\hhat_n}(t) - f(t))^2] \leq C(a) r_n(t,f). \]
\end{theorem}

\begin{remark}
  By specifying Theorem~\ref{THM:ADAPTATION} with the $\sinc$-kernel and $\hupp_n = 1$, one obtains an adaptive estimator attaining the optimal rate of convergence over functions bounded by $M$ in a Sobolev ellipsoid up to an extra logarithmic factor.
  A logarithmic loss for adaptation is commonly accepted and even known to be indispensable for pointwise estimation in the non-private framework \cite{brown1996constrained}.
\end{remark}
 \section{Discussion}

We have suggested an approach to adaptive kernel density estimation via Lepski's method in the framework of local approximate differential privacy.
Although we have studied its theoretical properties in the prototypical example of univariate density estimation only, our methodology should also be transferable to the multivariate case.
We also conjecture that it might be possible to extend our results to the case of general linear functionals (different from pointwise evaluation of the density function at a fixed point) as investigated in~\cite{goldenshluger2000adaptive} via Lepski's method in a inverse problem setup.
Furthermore, our methodology might be applicable to obtain local private estimation procedures in functional data analysis.
However, a lot of questions remain open:
One drawback of our approach is that the perturbation by a Gaussian process provides only approximate differential privacy and cannot be extended to pure differential privacy.
The creation of new methods for kernel estimators that overcome this restriction provides a further direction for future research.
Moreover, the optimal power of the logarithmic factor in the adaptive rate of convergence deserves further investigation as well as the behaviour of the minimax optimal rates in terms of the privacy parameters $\privpar$ and $\beta$.
 
\appendix

\section{Proofs of Section~\ref{SEC:PRIVACY}}

\subsection{Proof of Proposition~\ref{PROP:UNIV:LAPLACE:MECH}}

Let $A \in \Borel(\R)$ be arbitrary.
It has to be shown that
\[ \int_A  \frac{1}{2b} \exp\left( - \frac{\lvert z - g(x)\rvert}{b} \right) \dd z \leq e^\privpar \int_A  \frac{1}{2b} \exp\left( - \frac{\lvert z - g(x^\prime)\rvert}{b} \right) \dd z + \beta \]
for any $x,x^\prime \in \Xf$.
By the triangle inequality this holds true if
\[ \int_A  \frac{1}{2b} \exp\left( - \frac{\lvert z - g(x)\rvert}{b} \right) \dd z \leq e^{\privpar-\frac{\lvert g(x)-g(x^\prime) \rvert}{b}} \int_A  \frac{1}{2b} \exp\left( - \frac{\lvert z - g(x)\rvert}{b} \right) \dd z  + \beta, \]
and the latter holds true as soon as $1 \leq \exp \left( \privpar - \Delta(g)/b \right) + \beta$ which is equivalent to $b \geq \Delta(g)/(\privpar - \log(1- \beta))$.

\subsection{Proof of Proposition~\ref{PROP:PRIV:MULTI}}

We have to show that
\begin{equation*}
   \Pb^{Z|X=x}(A) \leq \exp(\privpar) \Pb^{Z|X=x^\prime}(A) + \beta
\end{equation*}
for all $A \in \Borel(\R^m)$.
This condition is satisfied if the set where the ratio $\dd \Pb^{Z|X=x}/\dd \Pb^{Z|X=x^\prime}$ exceeds $\exp(\privpar)$ has probability bounded by $\beta$ under $\Pb^{Z|X=x}$.
We have
\begin{align*}
 \frac{\dd \Pb^{Z|X=x}(z)}{\dd \Pb^{Z|X=x^\prime}(z)} = \exp \left( \frac{1}{2\sigma^2} \left[ -(z-g(x))^\transposed \Sigma^{-1}(z-g(x))  + (z-g(x^\prime))^\transposed \Sigma^{-1}(z-g(x^\prime)) \right] \right),
\end{align*}
and the condition $\frac{\dd \Pb^{Z|X=x}(z)}{\dd \Pb^{Z|X=x^\prime}(z)} > \exp(\privpar)$ is equivalent to
\begin{align*}
 \frac{1}{2\sigma^2} \left[ -(z-g(x))^\transposed \Sigma^{-1}(z-g(x))  + (z-g(x^\prime))^\transposed \Sigma^{-1}(z-g(x^\prime)) \right] > \privpar
\end{align*}
which in turn can be reformulated as
\begin{equation}\label{EQ:DEF:COND:OMEGA}
  2z^\transposed \Sigma^{-1} (g(x)-g(x^\prime)) + (g(x^\prime))^\transposed \Sigma^{-1} g(x^\prime) - g(x)^\transposed \Sigma^{-1} g(x) > 2\sigma^2 \privpar.
\end{equation}
Set $\Omega = \{ z \in \R^m : \eqref{EQ:DEF:COND:OMEGA} \text{ holds} \}$ and let $\xi$ denote a $\Nc(0,I_m)$ distributed random variable where $I_m$ denotes the $m \times m$-dimensional identity matrix.
Then
\begin{align*}
 \Pb^{Z|X=x}(\Omega) &= \Pb ( 2\sigma \xi^\transposed \Sigma^{-\frac{1}{2}} (g(x)-g(x^\prime)) > 2\sigma^2 \privpar - (g(x)-g(x^\prime))^\transposed \Sigma^{-1}(g(x)-g(x^\prime)) )\\
 &\leq \Pb ( \xi^\transposed \Sigma^{-1/2} (g(x)-g(x^\prime)) > \sigma \privpar - \frac{\Delta^2}{2\sigma} )\\
 &\leq \Pb ( \Delta \nu > \sigma \privpar - \frac{\Delta^2}{2\sigma} )\\
 &= \Pb ( \nu > \frac{\sigma \privpar}{\Delta} - \frac{\Delta}{2\sigma} )
\end{align*}
where $\nu$ is a univariate standard Gaussian random variable.
We now use the standard estimate $\Pb(\nu \geq t) \leq e^{-t^2/2}/2$ whose right-hand side is smaller than $\beta \in (0,1/2)$ if
$t^2 \geq -2 \log(2\beta)$.
We apply this estimate with $t=\frac{\sigma \privpar}{\Delta} - \frac{\Delta}{2\sigma}$, and thus $\Pb^{Z|X=x}(\Omega) \leq \beta$ if
\begin{align*}
	\frac{\sigma \alpha}{\Delta} - \frac{\Delta}{2\sigma} \geq \sqrt{2\log(1/(2\beta))},
\end{align*}
and this holds at least if
\[ \sigma \geq \frac{\Delta}{\privpar} \sqrt{2 \log(1/(2\beta)) + 2\privpar}.  \]

\subsection{Proof of Lemma~\ref{LEM:COMP}}

Let $A \in \Zs_1 \otimes \Zs_2$ be a measurable set.
Denote $A_{z_1} = \{ z_2 \in \Zf_2 : (z_1,z_2) \in A \}$ which is measurable.
By Cavalieri's principle and the independence assumption
\begin{align*}
 \Pb^{Z | X=x} (A) &= \int_{\Zf_1} \Pb^{Z_2|X=x}(A_{z_1}) \Pb^{Z_1|X=x}(\dd z_1)\\
 &\leq \int_{\Zf_1} (e^{\privpar_2} \Pb^{Z_2|X=x^\prime}(A_{z_1}) \wedge 1 + \beta_2) \Pb^{Z_1|X=x}(\dd z_1)\\
 &=  \int_{\Zf_1}  (e^{\privpar_2}\Pb^{Z_2|X=x^\prime}(A_{z_1}) \wedge 1) \Pb^{Z_1|X=x}(\dd z_1) + \int_{\Zf_1} \beta_2 \Pb^{Z_1|X=x}(\dd z_1).
\end{align*}
Now put $\Omega = \{ \dd \Pb^{Z_1|X=x}/\dd \Pb^{Z_1|X=x^\prime} \leq e^{\privpar_1} \} \subseteq \Zf_1$.
Then $\Pb^{Z_1|X=x}(\Omega^\complement) \leq \beta_1$ since otherwise there would be a contradiction to approximate differential privacy.
Hence,
\begin{align*}
 \Pb^{Z | X=x} (A) &\leq \int_{\Zf_1 \cap \Omega} e^{\privpar_1 + \privpar_2}\Pb^{Z_2|X=x^\prime}(A_{z_1})\Pb^{Z_1|X=x^\prime}(\dd z_1)  + \int_{\Zf_1 \cap \Omega^\complement} \Pb^{Z_1|X=x}(\dd z_1) + \beta_2\\
 &\leq e^{\privpar_1 + \privpar_2} \Pb^{Z | X=x^\prime} (A)  + \Pb^{Z_1|X=x}(\Omega^\complement) + \beta_2\\
 & \leq e^{\privpar_1 + \privpar_2} \Pb^{Z | X=x^\prime} (A)  + \beta_1 + \beta_2
\end{align*}
which shows the claim assertion.

\section{Proofs of Section~\ref{SEC:MINIMAX}}

\subsection{Proof of Proposition~\ref{PROP:UPPER}}

The bias-variance decomposition for the estimator $\fhat_h(t)$ is
  \[ \Eb [ (\fhat_h(t) - f(t))^2 ] = (f_h(t) - f(t))^2 + \Eb [ (\fhat_h(t) - f_h(t))^2 ] \]
  where $f_h(t) = \Eb [\fhat_h(t)]$.
  We begin with the analysis of the bias.
 First recall that
 \[ f(t) = \frac{1}{2\pi} \int e^{-\ii t\omega} \Fc[f](\omega) \dd \omega, \]
 and due to centredness of the error added by the privacy mechanism
 \begin{align*}f_h(t) &= \int \frac{1}{h} K_{\sinc} \left( \frac{u-t}{h} \right) f(u) \dd u\\
  &= \frac{1}{2\pi h} \int \Fc\left[ K_{\sinc}\left( \frac{\cdot - t}{h} \right) \right](\omega) \Fc[f](\omega) \dd \omega\\
  &= \frac{1}{2\pi} \int e^{-\ii t \omega}\Fc\left[ K_{\sinc} \right](h\omega) \Fc[f](\omega) \dd \omega.
 \end{align*}
 Thus, using that $\Fc\left[ K_{\sinc} \right](\cdot) = \1_{[-\pi,\pi]}(\cdot)$, we obtain
 \begin{align*}
   (f_h(t)-f(t))^2 &= \frac{1}{4\pi^2} \left( \int_{\R} e^{-\ii t \omega} \left[ 1 - \Fc\left[ K_{\sinc} \right](h\omega) \right] \Fc[f](\omega) \dd \omega \right)^2\\
   &= \frac{1}{4\pi^2} \left( \int_{\R} e^{-\ii t \omega} \1_{ \{ \lvert \omega \rvert > 1/h \} } \Fc[f](\omega) \dd \omega \right)^2\\
   &\leq \frac{1}{4\pi^2} \int \lvert \Fc[f](\omega) \rvert^2 \lvert \omega \rvert^{2s} \dd \omega \cdot \int \1_{ \{ \lvert \omega \rvert > 1/h \} } \lvert \omega \rvert^{-2s} \dd \omega\\
   &\leq \frac{2\pi L^2}{4\pi^2} \cdot \frac{2}{2s-1} h^{2s-1} = C(L,s)h^{2s-1}.
 \end{align*}
 Let us now consider the variance, where we have to distinguish between the case of Laplace mechanism and Gaussian mechanism.
 We denote
 \[ \ftilde_h(t) = \frac{1}{nh} \sum_{i=1}^n K_{\sinc}\left( \frac{X_i-t}{h} \right).\]
For the Laplace mechanism, we have by denoting $\xi \sim \Laplace(1)$ that
\begin{align*}
   \Eb [ (\fhat_h(t) - f_h(t))^2 ] &= \Var(\widetilde f_h) + \frac 1 n \Var \left( \frac{\Delta(K_{\sinc}((t-\cdot)/h)/h)}{\privpar - \log(1-\beta)} \, \xi\right)\\
   &\leq \frac{\lVert f \rVert_\infty \int K_{\sinc}^2(u) \dd u}{nh} + \frac{ 8\lVert K_{\sinc} \rVert^2_\infty}{nh^2(\privpar - \log(1-\beta))^2}\\
   &\leq C(\lVert f \rVert_\infty, K_{\sinc}, \privpar, \beta) \left[ \frac{1}{nh} + \frac{1}{nh^2} \right].
\end{align*}
In a similar fashion, for the Gaussian mechanism, now letting $\xi \sim \Nc(0,1)$, we have
 \begin{align*}
    \Eb [ (\fhat_h(t) - f_h(t))^2 ] &= \Var(\widetilde f_h) + \frac 1 n \Var \left( \frac{2  \norm{K_{\sinc}}_\infty \sqrt{2\log (1/2\beta)+2\alpha}} {\privpar h} \xi\right)\\
    &\leq \frac{\lVert f \rVert_\infty \int K_{\sinc}^2(u) \dd u}{nh} + \frac{4\lVert K_{\sinc}\rVert_\infty^2 (2\log(1/(2\beta)) + 2\privpar)}{nh^2\privpar^2}\\
    &\leq C(\lVert f \rVert_\infty, K_{\sinc}, \privpar, \beta)  \left[ \frac{1}{nh} + \frac{1}{nh^2} \right].
 \end{align*}
 The statement of the proposition follows now by combining the obtained bounds for squared bias and variance.

\subsection{Proof of Proposition~\ref{PROP:LOWER}}

Let $\ellhat$, $Q \in \Qc_\alpha$ be arbitrary as in the statement of the proposition.
  Define $\psi_n >0$ via $\psi_n^2=n^{-\frac{2s-1}{2s+1}}$.
	Let $f_{0,n}$, $f_{1,n}$ be two functions in $\Sc(s,L)$ (to be specified later on) such that $(f_{0,n}(t) - f_{1,n}(t))^2 \gtrsim \psi_n^2$.
	Using a general reduction argument (see \cite{tsybakov2009introduction}, Section~2.2) it can be shown that
	\begin{align*}
	     \sup_{f \in \Sc(s,L)} \psi_n^{-2} \Eb [(\ellhat - f(t))^2] &\geq \psi_n^{-2} \sup_{\theta \in \{0,1\} } \Eb [(\ellhat - f_{\theta, n}(t))^2]\\
	     &\gtrsim \psi_n^{-2} \inf_{\tau} \max_{\theta \in \{0,1\} } \Pb_\theta(\tau = 1 - \theta)
	\end{align*}
	where the infimum is taken over all $\{0,1\}$-valued test functions $\tau$ based on the observations $Z_1,\ldots,Z_n$ and $\Pb_\theta$ denotes the distribution of $Z_1,\ldots,Z_n$ if the true density of $X_1,\ldots,X_n$ is $f_{\theta,n}$.
In view of \cite{tsybakov2009introduction}, Theorem~2.2, Statement~(iii), the claim assertion follows if we can choose the functions $f_{0,n}$ and $f_{1,n}$ such that
\begin{enumerate}[(1),itemsep=8pt,leftmargin=3em]
  \item\label{LB:I:1} $f_{0,n}, f_{1,n} \in \Sc(s,L)$,
  \item\label{LB:I:2} $(f_{0,n}(t) - f_{1,n}(t))^2 \gtrsim \psi_n^2$, \quad and
  \item\label{LB:I:3} $\KL(\Pb_{0}, \Pb_{1}) \leq C < \infty$ for some $C$ independent of $n$.
\end{enumerate}
To construct such $f_{0,n}, f_{1,n}$ we use ideas from Section~6 of \cite{butucea2001exact} and refer to this paper also for some of the computations.
First, take a strictly positive probability density $f$ on $\R$ that is infinitely often continously differentiable.
Setting $\lVert f^{(s)} \rVert_2^2 = \frac{1}{2\pi} \int_\R \lvert \Fc[f](\omega) \rvert^2 \lvert \omega \rvert^{2s} \dd \omega$, we can further assume that $\lVert f^{(s)}\rVert_2 \leq L$.
Then, for $\delta \in (0,1/2)$, define the function $f_{0,n}$ by
\[ f_{0,n}(x) = f_0(x) = \left( \frac{\delta}{2} \right)^{\frac{1}{s + 1/2}} f \left( x \left( \frac{\delta}{2} \right)^{\frac{1}{s + 1/2}} \right). \]
In order to define the second hypothesis $f_{1,n}$ we consider the auxiliary function $\Ktilde_s$ as introduced on p.~26 of~\cite{butucea2001exact} (its construction in that paper is borrowed from \cite{tsybakov1998pointwise}).
In particular, note that $\Ktilde_s$ is compactly supported and satisfies $\lVert \Ktilde_s^{(s)} \rVert_2 \leq 1-\delta/2$ (thus $\Ktilde_s \in \Sc(s,1)$) and $\Ktilde_s(0) \geq (1-\delta)C(s) > 0$.
Set $h_n=(n(\exp(\privpar) - 1)^2)^{-1/(2s+1)}$, and put
\[ g_{n,s}(x) = cL h_n^{s-\frac{1}{2}} \Ktilde \left( \frac{x-t}{h_n} \right) \]
for some constant $c>0$.
Defining $\gamma_{n,s} = \int g_{n,s}(x) \dd x < \infty$, set
\[ f_{n,1}(x) = f_{n,0}(x)(1-\gamma_{n,s}) + g_{n,s}(x). \]
We now check conditions \ref{LB:I:1}--\ref{LB:I:3} from above.

\smallskip

\paragraph{Verification of \ref{LB:I:1}:} The proof follows step by step along the lines of the one in \cite{butucea2001exact} and we omit the details.
We only record the fact that
\[ \gamma_{n,s} = cLh_n^{s+\frac{1}{2}} \int \Ktilde_s(u)\dd u = O(h_n^{s+ \frac{1}{2}}) \]
which will be used below.

\paragraph{Verification of \ref{LB:I:2}:}
We have
\begin{align*}
(f_{n,0}(t) - f_{n,1}(t))^2 &= (f_{n,0}(t) - (1-\gamma_n)f_{0,n}(t) - g(t))^2\\
&= \lvert \gamma_{n,s} f_{0,n}(t) - g(t) \rvert^2\\
&\geq \lvert \lvert g(t) \rvert - \lvert \gamma_{n,s} f_{0,n}(t)\rvert \rvert^2.
\end{align*}
Now, since $g_{n}(t) = C h_n^{s-\frac{1}{2}}$ and $\gamma_n = O(h_n^{s+\frac{1}{2}})$, the last expression inside the outer absolute values is greater than $Ch_n^{s-\frac{1}{2}}$ for sufficiently large $n$, say $n \geq n_0$.
Hence for $n \geq n_0$,
\[ (f_{n,0}(t) - f_{n,1}(t))^2 \geq C h_n^{2s-1} = C(\alpha) n^{-\frac{2s-1}{2s+1}} \]
which is the desired bound.

\paragraph{Verification of \ref{LB:I:3}:}
By Equation~(14) in \cite{duchi2018minimax} we have
\begin{equation}\label{EQ:KL:TV}
  \KL(\Pb_0, \Pb_1) \leq 4n (\exp(\privpar) - 1)^2 \TV^2(\Pb_0^{X_1}, \Pb_1^{X_1}).
\end{equation}
Now
\begin{align*}
\TV(\Pb_0^{X_1},\Pb_1^{X_1}) &= \int \lvert f_{n,0}(x) - f_{n,1}(x) \rvert \dd x\\
&= \int_\R \lvert - \gamma_{n,s} f_{n,0}(x) + g_n(x) \rvert \dd x\\
&\leq \gamma_{n,s} \int \lvert f_{n,0}(x) \rvert \dd x +  \int_\R \lvert g_n(x) \rvert \dd x\\
&\leq O(h_n^{s+\frac{1}{2}}) + C h_n^{s-\frac{1}{2}} \int \Ktilde \left( \frac{x-t}{h_n} \right) \dd x\\
&\leq O(h_n^{s+\frac{1}{2}}) + C (n(\exp(\privpar)-1)^2)^{-\frac{1}{2}} \int \Ktilde (u) \dd u\\
&\leq C (n(\exp(\privpar)-1)^2)^{-\frac{1}{2}}.
\end{align*}
for $n$ sufficiently large.
Thus, by \eqref{EQ:KL:TV}, for $n$ sufficiently large
\[ \KL(\Pb_0,\Pb_1) \leq Cn(\exp(\privpar) - 1)^2 \TV^2(\Pb_0^{X_1},\Pb_1^{X_1}) \leq C.  \]
 \section{Proofs of Section~\ref{SEC:ADAPTATION}}

\subsection{Proof of Proposition~\ref{PROP:UPPER:ORACLE}}

Under Assumption~\eqref{EQ:ASS:BOCHNER}, we have that $\sup_{0 < \eta \leq h} \lvert f_\eta(t) - f(t)\rvert^2$ converges to zero as $h \to 0$.
Let $n \geq 3$.
By definition of $v^2(\cdot)$, $\lambda(\cdot)$ and $\hlow_n = \log(\sqrt n)/\sqrt n$ (since $\hupp_n = 1$),
\[ v^2(\hlow_n)\lambda^2(\hlow_n) \geq \frac{M \int K^2(u)\dd u}{\sqrt n\log(\sqrt n)} + \frac{C_{\ppp\betaprime}^2}{\log(\sqrt n)} \cdot \kappa \log(\sqrt n/\log(\sqrt n)), \]
hence $\liminf_{n \to \infty} v(\hlow_n)\lambda(\hlow_n) > 0$, and the set in the definition of $\hast_n$ is non-empty provided that $n$ is sufficiently large.
Now, the bias-variance decomposition of the pseudo estimator is
\begin{align*}
    \Eb [(\fhat_{\hast_n}(t) - f(t))^2] &= (f_{\hast_n}(t) - f(t))^2 + \Var_f (\fhat_{\hast_n})\\
    &\leq (f_{\hast_n}(t) - f(t))^2 +  v^2(\hast_n)\\
    &\leq \frac{v^2(\hast_n)\lambda^2(\hast_n)}{4} + v^2(\hast_n)\\
    &\leq \frac{5}{4} v^2(\hast_n) \lambda^2(\hast_n).
\end{align*}

Let now $h_0$ be the minimizer in the definition of $r_n(t,f)$.
We distinguish the cases $h_0 < a \hast_n$ and $h_0 \geq a \hast_n$.
First, if $h_0 < a \hast_n$, then
\begin{align*}
  r_n(t,f) &= \sup_{0 \leq \eta \leq h_0} (f_\eta(x) - f(x))^2 + \frac{M \int K^2(u) \dd u\log (n)}{nh_0} + \frac{C_{\ppp\betaprime}^2 \log(n)}{nh_0^2}\\
  &\geq \frac{M \int K^2(u) \dd u\log (n)}{nh_0} + \frac{C_{\ppp\betaprime}^2 \log(n)}{nh_0^2}\\
  &\geq \frac{M \int K^2(u) \dd u\log (n)}{an\hast_n} + \frac{C_{\ppp\betaprime}^2 \log(n)}{na^2 (\hast_n)^2}\\
  &\geq C(a,\kappa) v^2(\hast_n)\lambda^2(\hast_n).
\end{align*}
If $h_0 \geq a \hast_n$, then by the very definition of $\hast_n$ we obtain
\begin{equation*}
   r_n(t,f) \geq \sup_{0 \leq \eta \leq h_0} (f_\eta(t) - f(t))^2 \geq \sup_{0 \leq \eta \leq a\hast_n} (f_{\eta}(t) - f(t))^2 > \frac{v^2(a\hast_n)\lambda^2(a\hast_n)}{4},
\end{equation*}
and thus $r_n(t,f)  \gtrsim v^2(\hast_n)\lambda^2(\hast_n)$ also in this case.

\subsection{Proof of Theorem~\ref{THM:ADAPTATION}}

We consider the risk decomposition
\begin{equation*}
 \Eb[(\fhat_{\hhat_n}(t)-f(t))^2] = \Eb[(\fhat_{\hhat_n}(t)-f(t))^2 \1_{ \{ \hhat_n \geq \hast_n \} }] + \Eb[(\fhat_{\hhat_n}(t)-f(t))^2 \1_{ \{ \hhat_n < \hast_n \} }],
\end{equation*}
and study the two terms on the right-hand side separately.

\smallskip

\noindent\emph{Analysis of the first term (Case $\hhat_n \geq \hast_n$)}. Note that the quantities $v(\cdot),\lambda(\cdot)$ satisfy $v(h) \geq v(h^\prime)$ and $\lambda(h) \geq \lambda(h^\prime)$ for $h^\prime \geq h$.
Thus, using the inequality $(a+b)^2 \leq 2a^2 + 2b^2$, we have for $h \leq h^\prime$ that
\begin{align*}
 \psi(h^\prime, h) &= v(h^\prime)\lambda(h^\prime) + v(h^\prime, h) \lambda(h)\\
 &\leq v(h) \lambda(h) + 2 \sqrt 2 v(h) \lambda(h)\\
 &= (1+2\sqrt 2)v(h) \lambda(h).
\end{align*}
By the definition of $\psi$ and $\hhat_n$, we obtain
\begin{align*}
 \lvert \fhat_{\hhat_n}(t) - \fhat_{\hast_n}(t) \rvert \1_{ \{\hhat_n \geq \hast_n \} } &\leq \psi(\hhat_n, \hast_n)\\
 &\leq \sup \{ \psi(\eta,\hast_n) : \eta \in \Hc_n, \eta \geq \hast_n \}\\
 &\leq (1+ 2\sqrt 2)v(\hast_n) \lambda(\hast_n).
\end{align*}
Hence (recall that we denote $f_h(t) = \Eb[\fhat_h(t)]$),
\begin{align*}
 \Eb[(\fhat_{\hhat_n}(t) - f(t))^2 \1_{ \{ \hhat_n \geq \hast_n \} }]&\\
 &\hspace{-9em}\leq 2 \Eb[(\fhat_{\hhat_n}(t) - \fhat_{\hast_n}(t))^2 \1_{ \{ \hhat_n \geq \hast_n \} }] + 2 \Eb[(\fhat_{\hast_n}(t) - f(t))^2]\\
 &\hspace{-9em}= 2 \Eb[(\fhat_{\hhat_n}(t) - \fhat_{\hast_n}(t))^2 \1_{ \{ \hhat_n \geq \hast_n \} }] + 2 \Eb[(\fhat_{\hast_n}(t) - f_{\hast_n}(t))^2] + 2 (f_{\hast_n}(t) - f(t))^2\\
 &\hspace{-9em}\lesssim v^2(\hast_n)\lambda^2(\hast_n)
\end{align*}
where we used the bound $\Var(\fhat_{\hast_n}) \leq v^2(\hast_n)$ for the term $2 \Eb[(\fhat_{\hast_n}(t) - f_{\hast_n}(t))^2]$ and the definition of $\hast_n$ to bound the term $(f_{\hast_n}(t) - f(t))^2$.

\smallskip

\noindent\emph{Analysis of the second term (Case $\hhat_n < \hast_n$)}.
For $h, \eta \in \Hc_n$ with $\eta < h$, set
\[ B_n(t,h,\eta) = \{ \lvert \fhat_h(t) - \fhat_\eta(t) \rvert > \psi(h,\eta) \} .\]
Let $h$ in $\Hc_n$.
Then, by definition of $\hhat_n$,
\[ \{ \hhat_n = a^{-1}h \} \subseteq \bigcup_{\substack{\eta \in \Hc_n \\ \eta < h}} B_n(t,h,\eta), \]
and thus
\begin{align*}
    \{ \hhat_n < \hast_n \} &= \bigcup_{\substack{h \in \Hc_n\\h< \hast_n}} \{ \hhat_n = h \}\\
    &\subseteq \bigcup_{\substack{h \in \Hc_n\\h< a\hast_n}} \{ \hhat_n = a^{-1}h \}\\
    &= \bigcup_{\substack{h \in \Hc_n\\h< a\hast_n}} \bigcup_{\substack{\eta \in \Hc_n\\\eta < h}} B_n(t,h,\eta).
\end{align*}
We obtain
\begin{align*}
   \Eb [(\fhat_{\hhat_n}(t) - f(t))^2 \1_{\{ \hhat_n < \hast_n \} } ] &\leq \sum_{\substack{h \in \Hc_n \\ h < a\hast_n}} \Eb [(\fhat_{a^{-1}h}(t) - f(t))^2 \1_{ \{ \hhat_n = a^{-1}h \} } ]\\
   &\leq \sum_{\substack{h \in \Hc_n \\ h < a\hast_n}} \sum_{\substack{\eta \in \Hc_n\\ \eta < h}} \Eb [(\fhat_{a^{-1}h}(t) - f(t))^2 \1_{B_n(t,h,\eta)}].
\end{align*}
By definition of $\hast_n$, for all $\eta, h \in \Hc_n$ with $\eta < h \leq \hast_n$, it holds
\[ \lvert f_{\eta}(t) - f(t) \rvert \leq \frac{v(\hast_n)\lambda(\hast_n)}{2} \leq \frac{v(h)\lambda(h)}{2}. \]
Now, for $\eta < h \leq \hast_n$,
\begin{align*}
    B_n(t,h,\eta) &= \{ \lvert \fhat_h(t) - \fhat_\eta(t) \rvert > \psi(h,\eta) \}\\
    &= \{ \lvert \fhat_h(t) - \fhat_\eta(t) - (f_h(x) - f_\eta(t)) + f_h(t) - f_\eta(t) - f(t) + f(t) \rvert > \psi(h,\eta) \}\\
    &\subseteq \left\{ v(h)\lambda(h) + \left\lvert \frac{1}{n} \sum_{i=1}^n \zeta_i \right\rvert > \psi(h,\eta) \right\}\\
    &\subseteq \left\{ \left\lvert \frac{1}{n} \sum_{i=1}^n \zeta_i \right\rvert > v(h,\eta)\lambda(\eta) \right\}
\end{align*}
where $\zeta_i=\zeta_{i,h,\eta}= Z_{i,h}(t) - Z_{i,\eta}(t) - (f_h(t) - f_\eta(t))$.
Note that $\Eb \zeta_i=0$ and $\Var(\zeta_i) \leq nv^2(h,\eta)$.
Now, by the Cauchy-Schwarz inequality,
\begin{align*}
    \Eb[(\fhat_{\hhat_n}(t) - f(t))^2 \1_{ \{ \hhat_n < \hast_n \} }] &\leq \sum_{\substack{h \in \Hc_n\\h< a\hast_n }} \sum_{\substack{\eta \in \Hc_n\\\eta < h}} \Eb[ (\fhat_{a^{-1}h}(t) - f(t))^2 \1_{ \{ \lvert \frac{1}{n} \sum_{i=1}^n \zeta_i \rvert > v(h,\eta)\lambda(\eta) \}} ]\\
    &\hspace{-6em}\leq \sum_{\substack{h \in \Hc_n\\h< a\hast_n }} \sum_{\substack{\eta \in \Hc_n\\\eta < h}} \left( \Eb [ (\fhat_{a^{-1}h}(t) - f(t))^4] \right)^{1/2} \left( \Pb \left( \bigg\lvert \frac{1}{n} \sum_{i=1}^n \zeta_i \bigg\rvert > v(h,\eta)\lambda(\eta) \right) \right)^{1/2}.
\end{align*}
For the first term in the sum, we have
\begin{align*}
      \Eb [ (\fhat_{a^{-1}h}(t) - f(t))^4] &= \Eb [(\fhat_{a^{-1}h}(t) - f_{a^{-1}h}(t) + f_{a^{-1}h}(t) - f(t))^4]\\
      &\leq 8 \Eb [(\fhat_{a^{-1}h}(t) - f_{a^{-1}h}(t))^4] + 8(f_{a^{-1}h}(t) - f(t))^4.
\end{align*}
Putting $\zeta_i^\prime = Z_{i,a^{-1}h}(t) - f_{a^{-1}h}(t)$, we have
\begin{align*}
	\Eb \left[ (\fhat_{a^{-1}h}(t) - f_{a^{-1}h}(t))^4 \right] &= \Eb \left[ \left( \frac{1}{n} \sum_{i=1}^n \zeta_i^\prime \right)^4 \right]\leq \frac{\Eb [(\zeta_i^\prime)^4]}{n^3} + \frac{3(\Eb [(\zeta_i^\prime)^2])^2}{n^2}.
\end{align*}
On the one hand,
\begin{align*}
\Eb [(\zeta_i^\prime)^4] &\lesssim \frac{C_{\ppp\betaprime}^4}{a^{-4}h^4} + \frac{1}{a^{-4}h^4} \Eb \left[ \left( K\left( \frac{X-t}{a^{-1}h} \right) - \Eb \left[ K\left( \frac{X-t}{a^{-1}h} \right) \right] \right)^4  \right]\\
&\lesssim \frac{C_{\ppp\betaprime}^4}{a^{-4}h^4} + \frac{8}{a^{-4}h^4} \Eb \left[ \left( K\left( \frac{X-t}{a^{-1}h} \right) \right)^4 \right] + \frac{8}{a^{-4}h^4} \left( \Eb \left[ K\left( \frac{X-t}{a^{-1}h} \right) \right] \right)^4\\
&\lesssim \frac{1}{a^{-4}h^4} + \frac{1}{a^{-3}h^3} + 1,
\end{align*}
on the other hand
\begin{align*}
\Eb [(\zeta_i^\prime)^2] \lesssim \frac{C_{\ppp\betaprime}^2}{a^{-2}h^2} + \frac{1}{a^{-1}h}.
\end{align*}
Hence,
\[ \Eb [ (\fhat_{a^{-1}h}(t) - f_{a^{-1}h}(t))^4 ]  \leq C v^4(a^{-1}h). \]
Moreover, for $a^{-1}h < \hast_n$,
\begin{align*}
   (f_{a^{-1}h}(t) - f(t))^4 \leq  \frac{v^4(\hast_n)\lambda^4(\hast_n)}{16}
\end{align*}
by the very definition of $\hast_n$.
Thus, altogether,
\[ \Eb[ (\fhat_{a^{-1}h}(t) - f(t))^4] \leq C(v^4(a^{-1}h)  + v^4(\hast_n) \lambda^4(\hast_n) ), \]
and by the monotonicity of $v(\cdot)$ and $\lambda(\cdot)$, for $\eta < h \leq \hast_n$
\[ \Eb[ (\fhat_{a^{-1}h}(t) - f(t))^4] \leq C \lambda^4(\eta)v^4(a^{-1}h). \]
Write $\zeta_i = \zeta_{i}^{(1)} + \zeta_{i}^{(2)}$
where $\zeta_{i}^{(1)} = K_h(X_i - t) - K_\eta(X_i - t) - (f_h(t) - f_\eta(t))$ and $\zeta_{i}^{(2)} = \frac{C_{\ppp \betaprime}}{\sqrt 2 h} \xi_{i,h} + \frac{C_{\ppp \betaprime}}{\sqrt 2 \eta} \xi_{i,\eta}$ with $\xi_{i,h}, \xi_{i,\eta}$ i.i.d.\,$\sim \Laplace(1)$ or $\zeta_{i}^{(2)} = \frac{C_{\ppp \betaprime} }{h} \xi_{i,h} + \frac{C_{\ppp \betaprime} }{\eta} \xi_{i,\eta}$ with $\xi_{i,h}, \xi_{i,\eta}$ i.i.d.\,$\sim \Nc(0,1)$ for $i=1,\ldots,n$.
We have
\begin{align*}
  \Pb \left( \left\lvert \frac{1}{n}\sum_{i=1}^n \zeta_i \right\rvert > v(h,\eta)\lambda(\eta) \right) &\leq \Pb \left( \left\lvert \frac{1}{n}\sum_{i=1}^n \zeta_i^{(1)} \right\rvert > \frac{v(h,\eta)\lambda(\eta)}{2} \right)\\
  &\hspace{5em}+ \Pb \left( \left\lvert \frac{1}{n}\sum_{i=1}^n \zeta_i^{(2)} \right\rvert > \frac{v(h,\eta)\lambda(\eta)}{2} \right).
\end{align*}
Consider $\Pb \left( \left\lvert \frac{1}{n}\sum_{i=1}^n \zeta_i^{(1)} \right\rvert > \frac{v(h,\eta)\lambda(\eta)}{2} \right)$ first.
By Bernstein's inequality (see Lemma~\ref{LEM:BERNSTEIN}) with $b=4\lVert K \rVert_\infty/\eta$,
\begin{align*}
  \Pb &\left( \left\lvert \frac{1}{n}\sum_{i=1}^n \zeta_i^{(1)} \right\rvert > \frac{v(h,\eta)\lambda(\eta)}{2} \right)\\
  &\leq 2\max \left\{  \exp \left( -\frac{nv^2(h,\eta)\lambda^2(\eta)}{4 nv^2(h,\eta)} \right), \exp \left( -\frac{n v(h,\eta)\lambda(\eta) \eta}{32\lVert K \rVert_\infty} \right) \right\}\\
  &= 2 \max \left\{  \exp \left( -\frac{\lambda^2(\eta)}{4} \right), \exp \left( -\frac{n v(h,\eta)\lambda(\eta) \eta}{32\lVert K \rVert_\infty} \right) \right\}.
\end{align*}
Note that
\[ v(h,\eta) \geq \frac{C_{\ppp\betaprime}}{\sqrt n \eta}. \]
For any $h \in \Hc_n$ and $n$ large enough, it holds
\[ \sqrt n \geq \sqrt \kappa \log(\sqrt n) \geq \sqrt \kappa \log(\hupp_n \sqrt n) = \sqrt \kappa \log( \hupp_n/(1/\sqrt n))) \geq \sqrt \kappa \log(\hupp_n/h). \]
Thus
\begin{align}
    \Pb \left( \left\lvert \frac{1}{n}\sum_{i=1}^n \zeta_i^{(1)} \right\rvert > \frac{v(h,\eta)\lambda(\eta)}{2} \right) &\leq 2 \max \left\{  \exp \left( -\frac{\lambda^2(\eta)}{4} \right), \exp \left( - \frac{C_{\ppp\betaprime} \sqrt n \lambda(\eta)}{32 \lVert K \rVert_\infty} \right) \right\} \notag\\
&\leq 2 \exp \left( - \kappa \left( \frac{1}{4} \wedge \frac{ C_{\ppp\betaprime}}{32 \lVert K \rVert_\infty} \right) \log \left( \frac{\hupp_n}{\eta} \right) \right). \label{EQ:PROB:ZETA1}
\end{align}
For the probability in terms of $\zeta_i^{(2)}$, we consider now the Gaussian case first.
Using standard concentration results for the Gaussian distribution, we obtain
\[ \Pb \left( \left\lvert \frac{1}{n}\sum_{i=1}^n \zeta_i^{(2)} \right\rvert > \frac{v(h,\eta)\lambda(\eta)}{2} \right) \leq 2 \exp\left( - \frac{n^2 t^2}{2\sigma^2} \right) \]
where $t= v(h,\eta)\lambda(\eta)/2$ and $\sigma^2$ denotes the variance of the Gaussian random variable $\sum_{i=1}^n \zeta_i^{(2)}$.
Then,
\[ \frac{n^2 t^2}{2\sigma^2} = \frac{n^2v^2(h,\eta) \lambda^2(\eta)}{8 \sigma^2} \geq \frac{\lambda^2(\eta)}{8} \geq \frac{\kappa \log(\hupp_n/\eta)}{8}. \]
Thus,
\begin{equation}\label{EQ:PROB:ZETA2}
  \Pb \left( \left\lvert \frac{1}{n}\sum_{i=1}^n \zeta_i^{(2)} \right\rvert > \frac{v(h,\eta)\lambda(\eta)}{2} \right) \leq 2 \exp \left( - \frac{\kappa \log(\hupp_n/\eta)}{8} \right).
\end{equation}
Combining \eqref{EQ:PROB:ZETA1} and \eqref{EQ:PROB:ZETA2}, we obtain for the Gaussian case
\[ \Pb \left( \left\lvert \frac{1}{n}\sum_{i=1}^n \zeta_i \right\rvert > v(h,\eta)\lambda(\eta) \right) \leq 4\exp \left( - \left( \frac{\kappa}{8} \wedge \frac{\kappa C_{\ppp\betaprime}}{32 \lVert K \rVert_\infty} \right) \log(\hupp_n/\eta) \right), \]
and we denote $\kappap = \frac{\kappa}{8} \wedge \frac{\kappa C_{\ppp\betaprime}}{32 \lVert K \rVert_\infty}$.

Let us now consider the probability in terms of $\zeta_i^{(2)}$ for the Laplace case which is a little bit more involved since the sum of two Laplace random variables is not Laplace anymore.
We decompose
\begin{align*}
  \Pb \left( \left\lvert \frac{1}{n}\sum_{i=1}^n \zeta_i^{(2)} \right\rvert > \frac{v(h,\eta)\lambda(\eta)}{2} \right) &\leq \Pb \left( \left\lvert \frac{C_{\ppp\betaprime}}{\sqrt 2 nh}\sum_{i=1}^n \xi_{i,h} \right\rvert > \frac{v(h,\eta)\lambda(\eta)}{4} \right)\\
  &\hspace{1em}+ \Pb \left( \left\lvert \frac{C_{\ppp\betaprime}}{\sqrt 2 n\eta}\sum_{i=1}^n \xi_{i,\eta} \right\rvert > \frac{v(h,\eta)\lambda(\eta)}{4} \right).
\end{align*}
Consider only the first probability on the right-hand side, the bound for the second one following analogously.
By Bernstein's inequality (see Lemma~\ref{LEM:BERNSTEIN}, take the version with control on the moments applied with $t= v(h,\eta)\lambda(\eta)/4$, $v^2 = C_{\ppp\betaprime}^2/h^2$ and $b = C_{\ppp\betaprime}/h$)
\begin{align*}
  \Pb \left( \left\lvert \frac{C_{\ppp\betaprime}}{\sqrt 2 nh}\sum_{i=1}^n \xi_{i,h} \right\rvert > \frac{v(h,\eta)\lambda(\eta)}{4} \right) &\leq 2 \max \left\{ \exp \left( - \frac{nt^2}{4v^2} \right), \exp \left( - \frac{nt}{4b} \right) \right\}\\
  &\leq 2 \max \left\{ \exp \left( - \frac{\lambda^2(\eta)}{64} \right), \exp \left( - \frac{\sqrt n \lambda(\eta)}{16} \right)  \right\},
\end{align*}
and hence by using $\sqrt n \geq \sqrt \kappa \log(\hupp_n/\eta)$,
\begin{align*}
    \Pb \left( \left\lvert \frac{C_{\ppp\betaprime}}{\sqrt 2 nh}\sum_{i=1}^n \xi_{i,h} \right\rvert > \frac{v(h,\eta)\lambda(\eta)}{4} \right) \leq 2 \exp \left( - \frac{\kappa}{64} \log \left( \frac{\hupp_n}{\eta} \right) \right).
\end{align*}
Finally, we obtain with $\kappap = \frac{\kappa}{64} \wedge \frac{\kappa C_{\ppp\betaprime}}{32 \lVert K \rVert_\infty}$ that
\[ \Pb \left( \left\lvert \frac{1}{n}\sum_{i=1}^n \zeta_i^{(2)} \right\rvert > \frac{v(h,\eta)\lambda(\eta)}{2} \right) \leq 4 \exp \left( - \kappap \log \left( \frac{\hupp_n}{\eta} \right) \right) \]
in the Laplace case.
Note that
\[ \Pb \left( \left\lvert \frac{1}{n}\sum_{i=1}^n \zeta_i^{(2)} \right\rvert > \frac{v(h,\eta)\lambda(\eta)}{2} \right) \leq 4 \exp \left( - \kappap \log \left( \frac{\hupp_n}{\eta} \right) \right) \]
for both cases with different choices of $\kappap$.
Now,
\begin{equation*}
  \Eb[(\fhat_{\hhat_n}(t) - f(t))^2 \1_{ \{ \hhat_n < \hast_n \} }] \lesssim \sum_{\substack{h \in \Hc\\h<a\hast_n}} \sum_{\substack{\eta \in \Hc\\\eta < h}} \lambda^2(\eta) v^2(a^{-1}h) \exp \left( - \frac{\kappap}{2} \log \left( \frac{\hupp_n}{\eta} \right) \right).
\end{equation*}
For sufficiently small $\gamma > 0$\footnote{Our calculations show that $\gamma > 0$ has to satisfy also that $\kappap/2 - \gamma -2 >0$. Such a choice is possible whenever $\kappap/2 - 2 > 0$ which holds for $\kappa$ large enough.}, we have
\begin{align*}
  \sum_{\substack{\eta \in \Hc_n\\ \eta < h}} \lambda^2(\eta)  \exp \left( - \frac{\kappap}{2} \log \left( \frac{\hupp_n}{\eta} \right) \right) &\lesssim \left( \frac{h}{\hupp_n} \right)^{\kappap/2 - \gamma} \sum_{\substack{\eta \in \Hc_n\\ \eta < h}} \log \left( \frac{\hupp_n}{\eta} \right) \left( \frac{\eta}{\hupp_n} \right)^\gamma\\
  &\lesssim \left( \frac{h}{\hupp_n} \right)^{\kappap/2 - \gamma} \sum_{j=0}^\infty j a^{-\gamma j}  \log(a)\\
  &\lesssim \left( \frac{h}{\hupp_n} \right)^{\kappap/2 - \gamma}.
\end{align*}
Recall that $v^2(h) \asymp \frac{1}{nh} + \frac{1}{nh^2}$.
Thus,
\begin{align*}
\Eb[(\fhat_{\hhat_n}(t) - f(t))^2 \1_{ \{ \hhat_n < \hast_n \} }] &\lesssim \sum_{\substack{h \in \Hc_n\\h<a\hast_n}} \left( \frac{h}{\hupp_n} \right)^{\kappap/2 - \gamma} v^2(a^{-1}h)\\
&\lesssim \frac{\hupp_n}{n} \sum_{\substack{h \in \Hc\\h<a\hast_n}} \left( \frac{h}{\hupp_n} \right)^{\kappap/2 - \gamma-1} + \frac{\hupp_n^2}{n\alpha^2} \sum_{\substack{h \in \Hc\\h<a\hast_n}} \left( \frac{h}{\hupp_n} \right)^{\kappap/2 - \gamma-2}.
\end{align*}
The sums on the right-hand side converge and the bound for the case $\hhat_n < \hast_n$ is negligible with respect to the upper bound $v^2(\hast_n)\lambda^2(\hast_n)$.

\section{Bernstein inequality}

The following version of the Bernstein inequality is taken from \cite{comte2015estimation}.

\begin{lemma}\label{LEM:BERNSTEIN}
Let $X_1,\ldots,X_n$ be i.i.d.\,random variables and put $S_n=\sum_{i=1}^n (X_i - \Eb[X_i])$.
Then, for any $t > 0$,
\begin{align*}
  \Pb (\lvert S_n - \Eb[S_n] \rvert \geq nt) &\leq 2 \exp \left( - \frac{nt^2}{2v^2 + 2b\eta} \right)\\
  &\leq 2 \max \left\{ \exp \left( -\frac{nt^2}{4v^2} \right), \exp \left( -\frac{nt}{4b} \right)  \right \}
\end{align*}
where $\Var(X_1) \leq v^2$ and $\lvert X_1 \rvert \leq b$ (or $\Eb [\lvert X_i \rvert^m] \leq \frac{m!}{2} v^2 b^{m-2} \text{ for } m \geq 2$).
\end{lemma}
 
\printbibliography

\end{document}